\documentclass[12pt]{article}
\usepackage{amsfonts,amsmath,graphicx}
\usepackage{color}
\def \a {{\alpha}}

\def \l {{\lambda}}

\def \s {{\sigma}}

\newlength{\pecettawidth}
\def \beq {\begin{eqnarray}}
\def \eeq {\end{eqnarray}}
\def \beqn {\begin{eqnarray*}}
\def \eeqn {\end{eqnarray*}}

\newcommand{\halmos}{\rule{1ex}{1.4ex}}

\newcounter{for}[section]

\numberwithin{equation}{section}

\newtheorem{itlemma}{Lemma}[section]
\newtheorem{itproposition}[itlemma]{Proposition}
\newtheorem{theorem}[itlemma]{Theorem}
\newtheorem{itcorollary}[itlemma]{Corollary}
\newtheorem{itremark}[itlemma]{Remark}
\newtheorem{itremarks}[itlemma]{Remarks}
\newtheorem{itdefinition}[itlemma]{Definition}
\newtheorem{itexample}[itlemma]{Example}

\newenvironment{fact}{\begin{itfact}\rm}{\end{itfact}}
\newenvironment{claim}{\begin{itclaim}\rm}{\end{itclaim}}
\newenvironment{lemma}{\begin{itlemma}}{\end{itlemma}}
\newenvironment{remark}{\begin{itremark}\rm}{\end{itremark}}
\newenvironment{remarks}{\begin{itremarks} \rm}{\end{itremarks}}
\newenvironment{corollary}{\begin{itcorollary}}{\end{itcorollary}}
\newenvironment{proposition}{\begin{itproposition}}{\end{itproposition}}
\newenvironment{definition}{\begin{itdefinition}\rm}{\end{itdefinition}}
\newenvironment{example}{\begin{itexample}\rm}{\end{itexample}}
\newenvironment{proof}{\noindent {\em Proof}.\ \
}{\hspace*{\fill}$\halmos$\medskip}
\newcommand{\be}[1]{\addtocounter{for}{1} \begin{equation}\label{#1}}
\newcommand{\ee}{\end{equation}}
\newcommand{\bl}[1]{\begin{lemma}\label{#1}}
\newcommand{\br}[1]{\begin{remark}\label{#1}}
\newcommand{\brs}[1]{\begin{remarks}\label{#1}}
\newcommand{\bt}[1]{\begin{theorem}\label{#1}}
\newcommand{\bd}[1]{\begin{definition}\label{#1}}
\newcommand{\bp}[1]{\begin{proposition}\label{#1}}
\newcommand{\bc}[1]{\begin{corollary}\label{#1}}
\newcommand{\bfact}[1]{\begin{fact}\label{#1}}
\newcommand{\bex}[1]{\begin{example}\label{#1}}
\newcommand{\ec}{\end{corollary}}
\newcommand{\efact}{\end{fact}}
\newcommand{\eex}{\end{example}}
\newcommand{\el}{\end{lemma}}
\newcommand{\er}{\end{remark}}
\newcommand{\ers}{\end{remarks}}
\newcommand{\et}{\end{theorem}}
\newcommand{\ed}{\end{definition}}
\newcommand{\ep}{\end{proposition}}
\newcommand{\epr}{\end{proof}}
\newcommand{\bpr}{\begin{proof}}
\newcommand{\bcl}[1]{\begin{claim}\label{#1}}
\newcommand{\ecl}{\end{claim}}

\newcommand{\ecs}{\end{corollary}}
\newcommand{\eers}{\end{exercise}}
\newcommand{\eexs}{\end{example}}
\newcommand{\eems}{\end{example}}
\newcommand{\els}{\end{lemma}}
\newcommand{\eles}{\end{lemmaex}}
\newcommand{\ets}{\end{theorem}}
\newcommand{\eds}{\end{definition}}
\newcommand{\eps}{\end{proposition}}

\newcommand{\bi}{\begin{itemize}}
\newcommand{\ei}{\end{itemize}}
\newcommand{\ben}{\begin{enumerate}}
\newcommand{\een}{\end{enumerate}}

\def\vbar{\mathchoice{\vrule height6.3ptdepth-.5ptwidth.8pt\kern-.8pt}
   {\vrule height6.3ptdepth-.5ptwidth.8pt\kern-.8pt}
   {\vrule height4.1ptdepth-.35ptwidth.6pt\kern-.6pt}
   {\vrule height3.1ptdepth-.25ptwidth.5pt\kern-.5pt}}
\def\fudge{\mathchoice{}{}{\mkern.5mu}{\mkern.8mu}}
\def\bbc#1#2{{\rm \mkern#2mu\vbar\mkern-#2mu#1}}
\def\bbb#1{{\rm I\mkern-3.5mu #1}}
\def\bba#1#2{{\rm #1\mkern-#2mu\fudge #1}}
\def\bb#1{{\count4=`#1 \advance\count4by-64 \ifcase\count4\or\bba A{11.5}\or
   \bbb B\or\bbc C{5}\or\bbb D\or\bbb E\or\bbb F \or\bbc G{5}\or\bbb H\or
   \bbb I\or\bbc J{3}\or\bbb K\or\bbb L \or\bbb M\or\bbb N\or\bbc O{5} \or
   \bbb P\or\bbc Q{5}\or\bbb R\or\bbc S{4.2}\or\bba T{10.5}\or\bbc U{5}\or
   \bba V{12}\or\bba W{16.5}\or\bba X{11}\or\bba Y{11.7}\or\bba Z{7.5}\fi}}

\def \qed {{\hspace*{\fill}$\halmos$\medskip}}

\def \s {\sigma}

\def \A {{\mathcal{A}}}

\def \O {{\cal{O}}}
\def \H {{\cal{H}}}

\def \F {{\cal{F}}}

\def \kG {{\cal{G}}}
\def \C {{\cal{C}}}

\def\sqr#1#2{{\vcenter{\vbox{\hrule height .#2pt
                             \hbox{\vrule width .#2pt height#1pt \kern#1pt
                                   \vrule width .#2pt}
                             \hrule height .#2pt}}}}

\def\pmb#1{\setbox0=\hbox{#1}%
   \kern-.025em\copy0\kern-\wd0
   \kern.05em\copy0\kern-\wd0
   \kern-.025em\raise.0433em\box0 }
\def\sqr#1#2{{\vcenter{\vbox{\hrule height.#2pt
     \hbox{\vrule width.#2pt height#1pt \kern#1pt
   \vrule width.#2pt}\hrule height.#2pt}}}}

\def\N{{\mathbb N}}
\def\Z{{\mathbb Z}}
\def\R{{\mathbb R}}

\def\E{{\mathbb E}}

\def\l{\lambda}

\def\a{\alpha}

\def\reff#1{(\ref{#1})}
\def \ind {\hbox{1\hskip -3pt I}}
\def\VE{{{{\cal X}_\textrm{ve}}(\gamma,\beta)}}
\def\KVE{K_{\rm VE}}
\def\ER{{{\mathcal X}_\text{e}}}

\newcommand {\acc}[1] {\left\{ {#1} \right\}}

\begin{document}
\title{On Diffusion Limited Deposition.}
\author{
  \renewcommand{\thefootnote}{\arabic{footnote}}
  A.\ Asselah, E.N.M.\ Cirillo, B.\ Scoppola, and E.\ Scoppola}
\date{}

\pagestyle{myheadings}
\markboth
    {Diffusion Limited Deposition}
    {Diffusion Limited Deposition}

\maketitle

\begin{abstract}
We propose a simple model of columnar growth through
{\it diffusion limited aggregation} (DLA). Consider a graph
$G_N\times\N$, where the basis has $N$ vertices
$G_N:=\{1,\dots,N\}$, and two vertices $(x,h)$ and $(x',h')$
are adjacent if $|h-h'|\le 1$. Consider there a simple
random walk {\it coming from infinity}
which {\it deposits} on a growing cluster as follows:
the cluster is a collection of columns, and the height of
the column first hit by the walk immediately grows by one unit.
Thus, columns do not grow laterally. 

We prove that there is a critical time scale $N/\log(N)$ 
for the maximal height of the piles, i.e., there exist
constants $\alpha<\beta$ such that the maximal pile height 
at time $\alpha N/\log(N)$ is of order $\log(N)$, while at time 
$\beta N/\log(N)$ is larger than $N^\chi$.
This suggests that a \emph{monopolistic regime} starts at such a 
time and only the highest pile goes on growing. 
If we rather consider a walk whose height-component
goes down deterministically, the resulting
\emph{ballistic deposition} has maximal height of order $\log(N)$ at time 
$N$.

These two deposition models, diffusive and ballistic, are also compared 
with uniform random allocation and Polya's urn. 
\end{abstract}

\newpage
\section{Introduction}
\label{sec-intro}
\paragraph{Motivation.}
A celebrated model of deposition via diffusion is
proposed in the early 80's by Witten and Sanders \cite{witten-sanders}.
The aggregate, denoted $A(K)$, made of $K$ sites of $\Z^d$ is
built inductively as follows. Choose $A(1)=\{0\}$ and assume $A(K)$.
Let $\partial A(K)$ denote its outer boundary. Informally, launch
a simple random walk, $n\mapsto S(n)$, far away from the origin, 
and stop it when it reaches $\partial A(K)$, say on random site $Y$. 
We set $A(K+1)=A(K)\cup \{Y\}$. In other words, if
$\tau_{\partial A(K)}$ is the time at which the walk hits $\partial A(K)$,
then for $y\in \partial A(K)$,
\[
P\big(A(K+1)=A(K)\cup \{y\}\big| A(K)\big)=\lim_{\|x\|\to\infty}
P_x\big( S(\tau_{\partial A(K)})=y\big| \tau_{\partial A(K)}<\infty\big).
\]
Simulations show that the cluster looks like a ramified tree with
long branches. Heuristically, the origin of reinforcement is clear. 
Think of the walk in terms of its radial component, which performs
an almost symmetric one-dimensional walk, and its transverse component.
Either the random walk sticks soon after reaching the outer radius of the
cluster, and it has to settle on a tip, or it takes time before
settling and its radial component diffuses, and has more
chances to visit the extremal shells,
hence increasing the probability of attaching a tip rather than
an inside site. This explains reinforcement, 
but does not explain why this reinforcement is enough to
produce a ramified tree structure. It is clear also, at the
heuristic level, that we face two problems: controlling the
number of tips in the growing cluster, 
and controlling in a quantitative
way the reinforcement of these tips.

One natural way to measure the dimension of the cluster
is to find the scaling of
the radius of $A(K)$, and look for $\bar d$ such that
\be{def-dim}
\text{Radius}(A(K))\sim K^{1/\bar d}.
\ee
If $A(K)$ were a ball, then $\bar d=d$, and the conjecture is
that $\bar d< d$. Now, physicists have a much sharper conjecture
\be{conj-2}
\bar d_c=d-\frac{d-1}{d+1}.
\ee
In dimension 2, $\bar d_c=5/3$, and simulations give $\bar d=1.7$.

Kesten in \cite{kesten-1,kesten-2,kesten-3} considers the problem, and shows
that the arms of the cluster
are not too long. More precisely, his result reads
\be{kesten-lower}
\bar d\ge
 \left\{ \begin{array}{ll}
& 3/2 \mbox{ for } d=2 \, ,\qquad (\bar d_c=2-1/3) \\
& 2 \mbox{ for } d= 3 \, ,\qquad(\bar d_c=5/2) \\
& d/2 \mbox{ for } d\ge 3 \, ,\quad(\bar d_c\le d-3/5).
\end{array} \right.
\ee

By reversing time, (see \cite{lawler} and assume $d\ge 3$)
one writes the probability of adding $Y=y$ to the cluster as
\be{def-harmonic}
P\big(A(K+1)=A(K)\cup \{y\}\big| A(K)\big)=
\frac{P_y\big(\tau_{\partial A(K)}=\infty)}
{\sum_{z\in \partial A}P_z\big(\tau_{\partial A(K)}=\infty)}.
\ee
The difficulty is to have estimate on the escape probability when
the set $A$ is not a sphere, or some simple geometric shape.
Let us mention an interesting result 
about holes in the DLA cluster, where a hole is a
finite maximal connected subset of the complement of $A(K)$.
Erbez-Wagner \cite{wagner} shows that in dimension two, 
almost surely the number of holes tends to infinity with $K$.

Barlow, Pemantle and Perkins in \cite{barlow-p-p} study DLA on
a regular $d$-ary tree where the conductance between edges joining
generation $n$ and $n+1$ is $\alpha^{-n}$ for $\alpha<1$.  
These authors show that the infinite cluster has
a unique infinite line of descent.
Even though there is an explicit formula for the
harmonic measure, the proof that
$r(A(K))$ scales like $K$ with normal fluctuations is
non-trivial.

Benjamini and Yadin in \cite{benjamini-yadin} propose another
toy model for DLA. They consider a cylinder $G_N\times \N$,
where the graph $G_N$ has constant degree, $N$ vertices, and
is fast mixing: the mixing-time should be less than
$\log^{2-\epsilon}(|G_N|)$ for some positive $\epsilon$
(the class of $d$-regular random graphs works).
They show that if we send $H\times |G_N|$ simple walks
from infinity, then the height of the aggregate
is larger than $H\log(\log(|G_N|))$ for any $H$ and $N$ large enough.

There is a two-dimensional model, the Hastings-Levitov model,
which takes advantage of the conformal invariance of two-dimensional
brownian motion, and Riemman's mapping Theorem to map
the complement of the cluster into the complement of the unit disk,
and then attach on the unit circle a {\it stick} 
at a random uniform angle. Recently
Norris and Turner \cite{norris-turner} have studied very precisely
the limiting cluster obtained
by iteration of randomly rotated conformal mappings.

In a series of three recent papers, Amir, Angel, Benjamini and Kozma
\cite{amir-1,amir-2,amir-3}
study DLA on $\Z$ with long-range 
random walks. The cluster is
no longer connected, and they discover many phase transitions in the
growth rate of the cluster according to the tail decay of 
the increment of the walk. 

Our model is a further simplification of Benjamini and Yadin's model
\cite{benjamini-yadin} in two ways:
(i) no lateral hairs are produced, and (ii) the basis graph has
no geometry.
In our toy model of DLA, the radial component 
does a one-dimensional random walk, and the transverse component
samples uniformly the section of our graph.
Still we believe that our model is interesting,
and one can answer some of
the following questions in a quantitative way. 
\begin{itemize}
\item  What is the origin of reinforcement?
\item  What is the critical height to overcome ?
\item  What are the different regimes in the cluster's growth?
\end{itemize}
\paragraph{Models.}
We shall consider two deposition models, 
\emph{diffusive deposition} and \emph{ballistic deposition}. 

We start with defining diffusive deposition. 
Our graph is a half--cylinder $G_N\times\N$, where the basis
has $N$ vertices $G_N:=\{1,\dots,N\}$, 
and two vertices $(x,h)$ and $(x',h')$
are adjacent if $|h-h'|=
1$. The set $G_N\times \{0\}$ is called {\it the ground}. 

Let $n\mapsto A(n)$ be the evolution of random subsets of $G_N\times\N$
that we call {\it the cluster}.
The cluster is built inductively with $A(0)=G_N\times \{0\}$.
For an integer $k$, the cluster $A(k)$ is made of columns, that is, 
\be{def-A}
A(k)=\bigcup_{i=1}^N \{i\}\times \{0,\dots,\sigma_i(k)\}\quad
\text{with}\quad\sum_{i=1}^N \sigma_i(k)=k.
\ee
We shall write for simplicity
$A(k)=(\sigma_1(k),\dots,\sigma_N(k))$.

Assume that $A(k)$ is built. We consider a simple
random walk $n\mapsto S_n=(X_n,Z_n)$ on our graph. In other words,
\begin{enumerate}
\item\label{mod01}
$\{X_n\}$ an i.i.d. sequence uniformly distributed on $G_N$;
\item\label{mod02}
$\{Z_{n+1}-Z_n\}$ i.i.d.\ uniformly on $\{-1,1\}$;
\item\label{mod03}
the initial condition $Z_0$ is above the maximal 
height of the cluster $A(k)$. For defininetness we take
$Z_0=\max_i\sigma_i(k)+1$.
\end{enumerate}

The following rule of {\it aggregation}, or {\it deposition},
makes the cluster grow. The walk $S_n$, roams until 
it hits the cluster $A(k)$. 
Let $(X^*,Z^*)$ be the hitting site on $A(k)$, and
necessarily $0\le Z^*\le \sigma_{X^*}(k)$. We build $A(k+1)$ by
increasing the height of column $X^*$ by one unit. That is
\begin{displaymath}
\sigma_i(k+1)=\sigma_i(k)\;\;\textrm{ for any }i\neq X^*
\;\;\;\textrm{ and }\;\;\;
\sigma_{X^*}(k+1)=\sigma_{X^*}(k)+1
\;\;.
\end{displaymath}
We shall also say that the walk \emph{attaches} to the column, or pile, 
at $X^*$.
The walk with the aggregation rule is called an {\it explorer}.
We shall denote by $P$ the probability associated with this process. 

In diffusive deposition there are two relevant phenomena: one is diffusion, the
other is deposition which happens instantly and this explains
the name {\it diffusion limited deposition}.

Ballistic deposition is defined similarly, with the same notation, 
but with a totally asymmetric walk $\{Z_{n+1}-Z_n=-1\}$.
One could consider a continuum of biased models with a drift parameter.

\paragraph{Definitions and notation.} 
We use $\sigma,\eta$ to denote configurations, i.e., 
$\sigma,\eta\in\mathbb{N}^N$, $\sigma=(\sigma_1,\sigma_2,\dots,\sigma_N)$.
We also let $|\sigma|:=\sum_{i=1}^N\sigma_i$.
The symbol $\bar\sigma$ will denote the configuration obtained by ordering 
the components of $\sigma$ so that 
$\bar\sigma_1\ge\bar\sigma_2\ge\cdots\ge\bar\sigma_N$. 
We call $\O_N$ the set of ordered configurations 
$\eta\in \N^N$, namely, 
such that $\eta_1\ge\eta_2\ge\dots\ge\eta_N$.

Given a configuration $\sigma$, we denote by $\zeta(\sigma)$ 
the \emph{height occupation} of $\sigma$, i.e.,
\begin{equation}
\label{heightocc}
\zeta_j(\sigma)=\sum_{i=1}^N \ind_{\{\sigma_i\ge j\}}
\;\;.
\end{equation}
Note that $\sum_{j\ge 1}\zeta_j(\sigma)=|\sigma|$, and that
$\zeta(\sigma)=\zeta(\bar\sigma)$. Given 
two configurations $\sigma$ and $\eta$ such that $|\sigma|=|\eta|$, 
we say that $\sigma$ is more {\it monopolistic} than $\eta$,
writing $\sigma\succ\eta$, when
\be{def-order}
\forall k=1,\dots,N\qquad \sum_{i=1}^k \bar\sigma_i \ge \sum_{i=1}^k \bar
\eta_i.
\ee
Equivalently, one realizes $\bar\eta$ from $\bar\sigma$
by moving particles from the highest columns to the lowest ones.

Urn models are paradigms of reinforcement phenomena (see for instance
the survey \cite{pemantle}),
and our deposition models actually can be stochastically compared
with urns with $N$ colors. We briefly recall Polya's urn with $N$ colors:
starting with one ball of each color, at each unit time one
draws a ball and put it back in the urn with an 
additional ball of the same color. 
Calling $\eta_i$, with $i=1,\dots,N$,  
the number of added balls of color $i$ after $|\eta|$ draws, 
the probability of drawing a ball of color $i$ is
\be{def-polya}
q_i^P(\eta)=\frac{\eta_i+1}{\sum_{j=1}^N\eta_j+N}\;\;.
\ee
We consider also a generalized urn 
by replacing the r.h.s.\ in \reff{def-polya}
by
$f(\eta_i)/\sum_{j=1}^Nf(\eta_j)$ 
where
$f:\N\to \R^+$ is a function such that $f(0)=1$.
When $f(x)=x^2+1$, we call the model the {\it quadratic urn}.
When $f\equiv 1$, we call the model the {\it uniform random allocation} 
and denote by $q_i^U(\eta)=1/N$ 
the corresponding probability of drawing a ball of color $i$ 
at any time. Finally, we say that a process $t\mapsto \sigma(t)$ 
is more {\it monopolistic}
than process $t\mapsto \eta(t)$ if there is a coupling of the processes
such that for any $t>0$ we have $\sigma(t)\succ \eta(t)$, if this is
the case initially. 


\paragraph{Main Results.}
In this Section we collect our main results. 
The first Theorem gives an estimate of the number of explorers 
necessary, in diffusive deposition,
to form a cluster with at least one column proportional to a power of $N$.

\bt{main-theo}
Consider diffusive deposition.
There are constants $\alpha<\beta$, such that 
almost surely, when $N$ is large enough
\be{result-1}
\max_{i\in G_N} \sigma_i\big(\frac{\alpha N}{\log(N)}\big)\le 3 \log(N),
\ee
and there exists a positive constant $\chi$ such that 
\be{result-2}
\max_{i\in G_N} \sigma_i\big(\frac{\beta N}{\log(N)}\big)\ge N^\chi.
\ee
\et
In ballistic deposition, we prove that the growth of the height of the cluster 
is much slower. Indeed, it is unlikely that $N$ explorers
produce a column of height $\log(N)$.

\bt{main-theo-bal}
Consider ballistic deposition. 
There exists a positive constant $A$ such that
almost surely, when $N$ is large enough
\be{result-3}
\max_{i\in G_N} \sigma_i(N)\le A\log(N).
\ee
\et
\br{rem-walk}
For the {\it radial} component, we could have
chosen $\{Z_{n+1}-Z_n\}$ i.i.d.\ with some finite range law
without affecting our results.
\er
Theorem~\ref{main-theo} occurs in a regime where less than
$N$ explorers are thrown in the graph. We call this the
{\it early regime} which is to be thought of as the configurations
where an additional explorer has good chances to settle on the
ground $G_N\times\{0\}$. To motivate other results, let us
explain the different steps leading to a column of height
$N^\chi$. A first step is to reach a {\it subcritical} height
$\log(N)/\log(\log(N))$. We obtain that a large number of
columns reach this height by comparison with random allocation.
We obtain interesting comparison with other urns,
with the observation that ballistic deposition looks
like Polya's urn with $N$ colors, whereas diffusive deposition
looks like a quadratic urn (see below \reff{def-polya} for
the definition). Then, one of these
subcritical columns reaches the {\it critical height} $\log(N)$.
Since our estimate requires the configuration to stay
in the early regime, one has to
bound the number of {\it critical columns}. We show that
the number of critical columns is less than $N^{1-2\chi}$ for
some positive $\chi$, and this implies that the evolution 
remains in the early regime as long as the highest column
has not crossed $N^\chi$. We now can state our comparison result.
\bt{theo-comparison}
Both deposition models (diffusive and ballistic) are more
monopolistic than Polya's urn, which itself is more monopolistic
than random allocation.
\et
The following corollary is a side result interesting
on its own right which seems new, to the best of our knowledge.
\bp{theo-polya} Polya's urn with $N$ colors is {\it monotone}
with respect to the order $\succ$.
\ep
\noindent

\paragraph{Related Models.}
There are many models of cluster growth similar in definition
to DLA. They differ according to the law of $Y$, the site we add
on the boundary of the cluster $A$. This can also be expressed
according to the site, say $X$, from where the
random walks are launched and lead to different phenomenology.
\begin{itemize}
\item If $X=0$, we rather define a dual model of {\it erosion}.
The cluster represents the eroded materia, and $A(0)=\emptyset$.
Each new walk starts at 0, and settles on the first visited site
outside the cluster (a site which we interpreted as being eroded). 
This is internal DLA, and was introduced by
Meakin and Deutch in\cite{meakin-deutch}. The cluster is spherical
as was first seen Lawler, Bramson and Griffeath in \cite{lawler-92}. 
The fluctuations were
studied in \cite{AG1,AG2,AG3} and independently in \cite{JLS1,JLS2,JLS3}.
\item If $X$ is uniformly drawn in the cluster, then 
Benjamini, Duminil-Copin, Kozma, Lucas in \cite{BDKL}
show that the cluster is spherical.
\item If $Y$ is uniform on the boundary of the cluster, then
this is the celebrated Eden model \cite{eden}, 
which was proposed in the '60, and studied first by Richardson
\cite{richardson}.
\item If particles do not erode immediately the materia, but do it
with an exponential clock, and if they can be activated again when
another walk stands on their site, this is Activated Random Walks.
This model has been introduced by Spitzer in the 70, and
much discussed in the physics literature as an example of 
{\it self-organized criticality}. This has been studied mathematically
by Rolla and Sidoravicius \cite{rolla-s} (and references therein), 
and recently by Sidoravicius and Teixera \cite{s-teixera} among others.
Recent efforts have focused on the case of an initial
condition drawn from a product Poisson measure. As one tunes the
density there is phase transition between settlement of explorers
(in any finite box), and their perpetual activity.
\end{itemize}

\paragraph{Pictures and simulations.} 
In order to illustrate our main results, we show some numerics. 
In particular, we emphasize the freezing 
phenomenon which leads to the monopolistic regime: 
after a given time the highets pile grows linearly catching 
all particles. We stress that simulations do not capture 
quantitative aspects of the problem (scalings or exponents), 
but serve merely as qualitative illustrations. 

For both the diffusive and the ballistic model 
we have simulated the systems for 
$N=50,100,200,\dots,1000$.
For the diffusive model we have considered also the cases 
$N=2000,4000,5000,6000,8000,10000$.
In all the cases averages have been computed over $10^4$ 
independent realizations of the process. 
We have checked in all the cases that the sample is large 
enought to get stable averages. 

\begin{figure}[h]
\begin{center}
\begin{picture}(200,120)(125,0)
\put(0,0)
{
\includegraphics[width=12.cm]{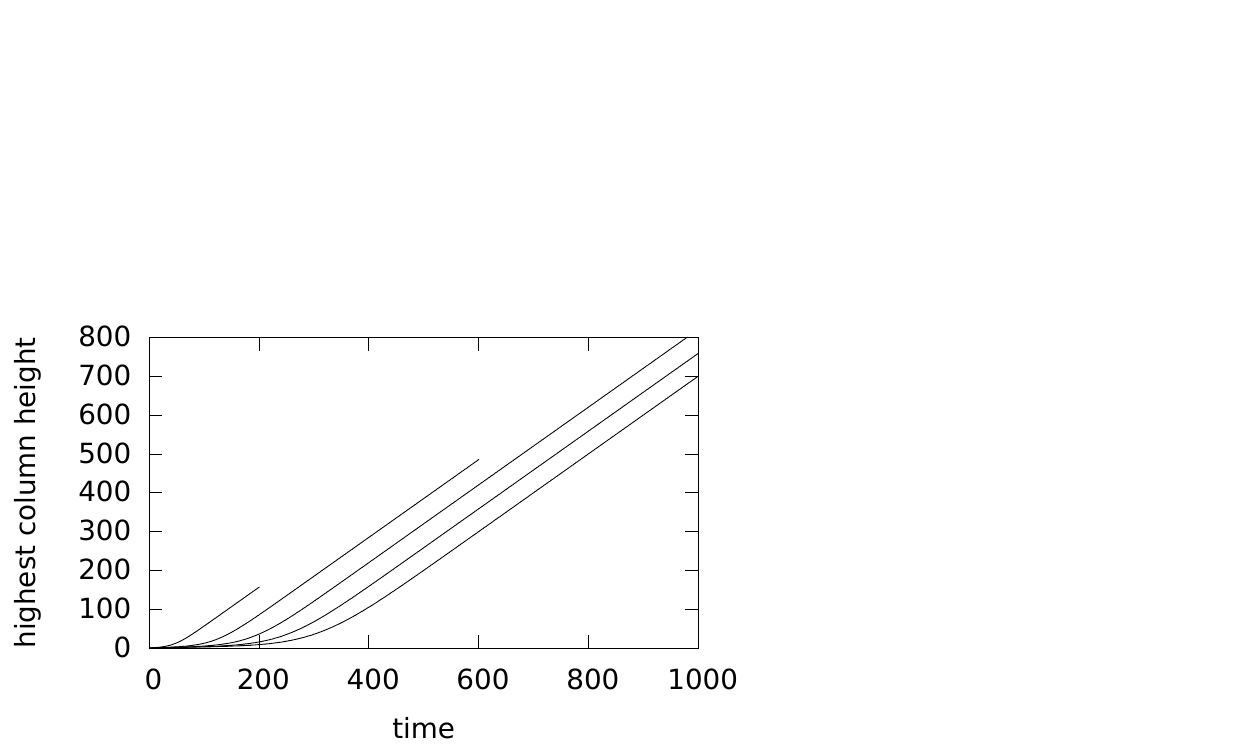}
}
\put(240,0)
{
\includegraphics[width=12.cm]{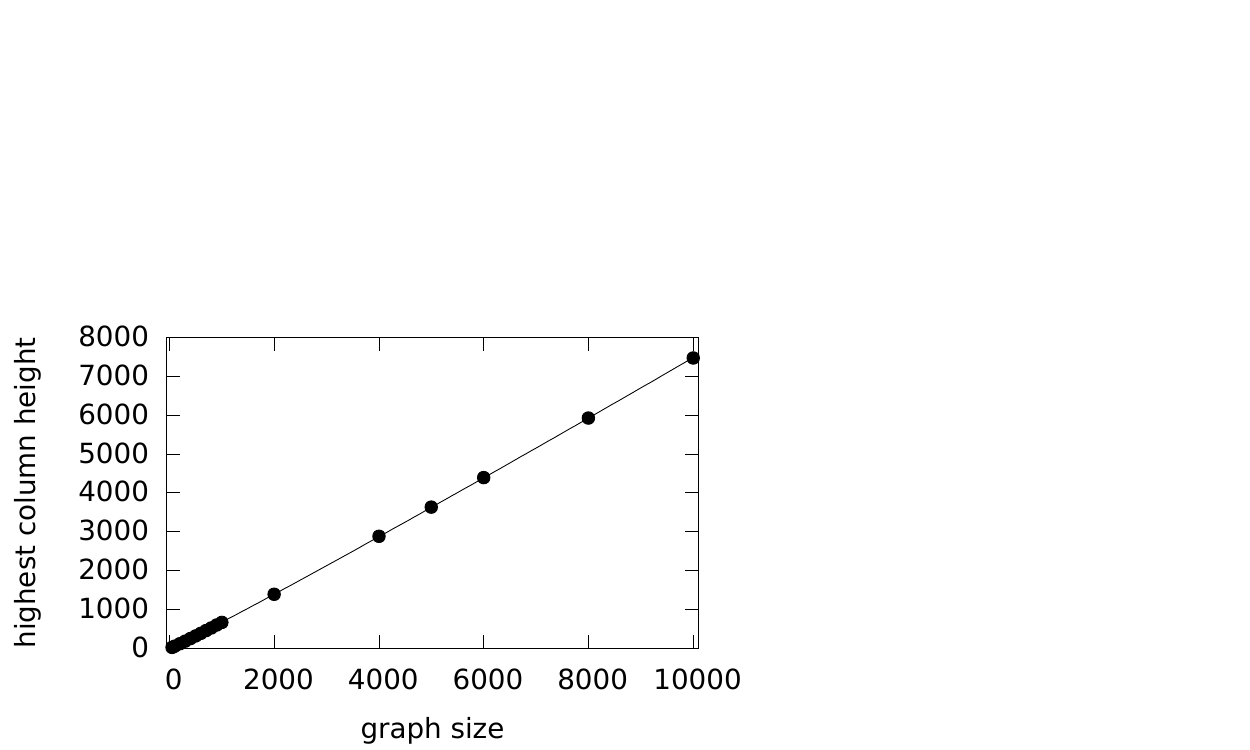}
}
\end{picture}
\caption{Numerical simulations for the diffusive model.
Left panel: the highest column height is plotted as function of time (number 
of explorers). The five plotted curves, from the left to the right, 
refer to $N=100,300,500,700,900$, respectively. 
Right panel: solid disks refer to the simulated 
highest column height at time $N$, namely, after $N$ 
explorers have been sent, for different values of the size of the graph $N$. 
The solid line is an eye--guide obtained by plotting the 
fitting function $0.498\times N^{1.044}$. 
}
\label{f:dif}
\end{center}
\end{figure}

\begin{figure}[h]
\begin{center}
\begin{picture}(200,120)(125,0)
\put(0,0)
{
\includegraphics[width=12.cm]{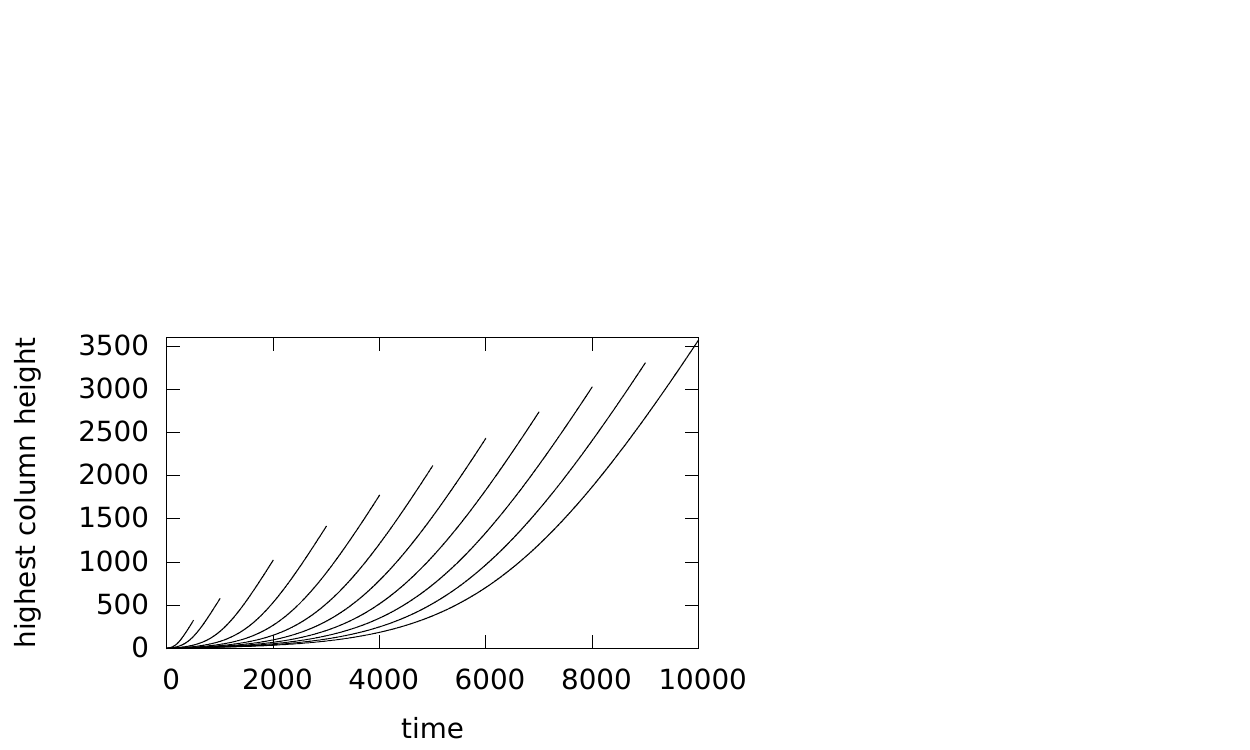}
}
\put(240,0)
{
\includegraphics[width=12.cm]{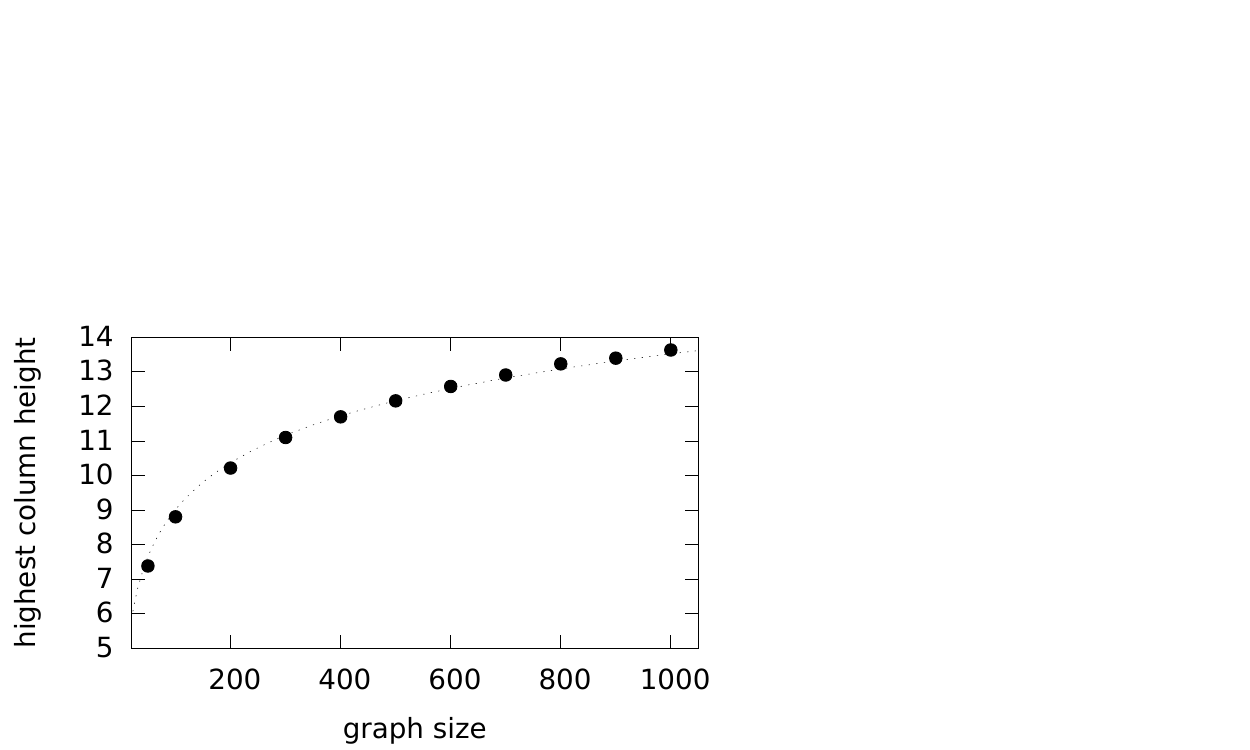}
}
\end{picture}
\caption{Numerical simulations for the ballistic model
for $N=50,100,200,\dots,1000$ with averages computed over $10^4$ 
independent realizations of the process. 
Left panel: the highest column height is plotted as function of time (number 
of explorers). The plotted curves, from the left to the right, 
refer to $N=50,100,200,\dots,1000$, respectively. 
Right panel: solid disks refer to the simulated 
highest column height at time $N$, namely, after $N$ 
explorers have been sent, for different values of the size of the graph $N$. 
The solid line is an eye--guide obtained by plotting the 
function $1.957\log(N)$. 
}
\label{f:bal}
\end{center}
\end{figure}

Simulations show clearly that the ballistic model 
reaches the monopolisitc regime much later than the 
diffusive one. Indeed in both models, compare the left panels in 
Figures~\ref{f:dif} and \ref{f:bal}, the height of the highest pile 
attains a linear behavior after an initial transient.  
This late time regime is the one in which all the particles are 
caught by the highest pile. Data show that the time length of the 
transient is much smaller in the diffusive model. 

We have also tested numerically our main results in 
Theorems~\ref{main-theo} and \ref{main-theo-bal}.
Indeed, we have computed, by averaging over different realizations 
of the process, the typical height of the highest pile at time 
$N$. Rigorous results suggest that this quantity should scale 
as a power law in the diffusive case and logarithmically in the 
ballistic one. The related numerical results are shown in 
the right panels in Figures~\ref{f:dif} and \ref{f:bal}. 
The qualitative agreement between simulations 
and theoretical results is striking. 
We stress again that the numerical results cannot be interpreted as a 
quantitave description of the model behavior since, for instance, 
too small values of the graph size $N$ have been considered. 

\begin{figure}[h]
\begin{center}
\begin{picture}(200,120)(0,0)
\put(-20,0)
{
\includegraphics[width=11.cm]{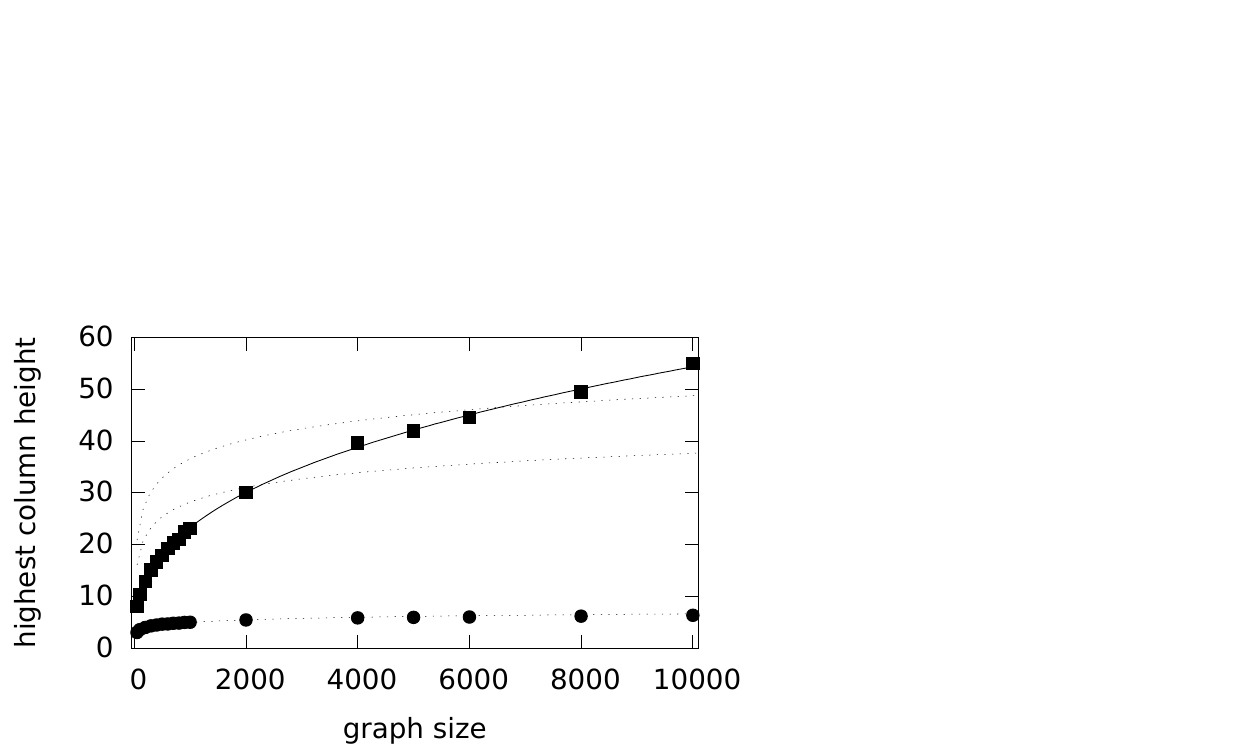}
}
\end{picture}
\caption{Numerical simulations for the diffusive model.
Solid disks and squares refer, respectively, to the simulated 
highest column height at times $N/\log(N)$ and 
$2 N/\log(N)$
for different values of the size of the graph $N$. 
The solid line is an eye--guide obtained by plotting the 
fitting function $1.847\times N^{0.367}$. The two 
dotted lines are the graph of the two functions 
$0.724\log(N)$, $4.085\log(N)$, and $5.292\log(N)$. 
}
\label{f:critico}
\end{center}
\end{figure}

Finally, in the diffusive case we have also tested numerically 
our results about the critical character of the 
time scale $N/\log(N)$. 
In Figure~\ref{f:critico} we compare the highest column 
height measured at times $N/\log(N)$ and $2N/\log(N)$.
In the first case the numerical (solid circles) data can be
perfectly fitted by a logarithmic function. In the latter case, 
on the other hand, the poor logarithmic fitting is opposed to a perfect 
power law one of the numerical data (solid squares). 
This result is in perfect agreement with the one 
proved in Theorem~\ref{main-theo} and, in particular, 
it suggests that $1<\alpha<\beta<2$. 
We stress that our numerics cannot in any case be considered 
quantitative, indeed, we have no clue to state that, by considering 
larger sizes of the graph, our numerical results would be confirmed. 

\paragraph{Plan.} The rest of the paper is organized as follows.
In Section~\ref{sec-tools} we present our main tool, which
is the probability an explorer hits the ground, as well
as a heuristic explanation for the logarithmic
scale of the {\it critical} height. In Section~\ref{sec-urns},
we present comparison with urns models. 
The proof of Proposition~\ref{theo-polya} 
is given in Section~\ref{s:com-pr}.
The proof of the Theorem~\ref{theo-comparison} 
is given in Section~\ref{s:com-dp}. In Section~\ref{sec-unit},
we establish our main tool, and related estimates. We study the
very early regime in Section~\ref{sec-VE}. Then, we study
the growth of cluster in Section~\ref{sec-column}, and the
reason why a large number of columns cannot overcome
height $\log(N)$. Finally, we gather all the needed estimates
to prove Theorem \ref{main-theo} in Section~\ref{sec-proof}.

\section{Key Tools and Sketch.}
\label{sec-tools}
The key to Theorem~\ref{main-theo} is an estimate of the
probability of attaching to a given column. 
Before, we need a lower bound on the probability 
of hitting first the ground.
\bl{lem-DD}
Consider diffusive deposition with a configuration $\sigma$
such that $|\sigma|<N/2$ explorers,
\be{DD-partition}
P_g(\sigma):=P\big(\text{ Explorer hits the ground }\big|\ \sigma\big)
\ge \exp\Big(-\frac{1}{N} \big(\sum_{j=1}^N \sigma_j(\sigma_j+1)\big)
\big(1+O(\frac{|\sigma|}{N}\big)\Big).
\ee
\el
The time spent on the slab $G_N\times \{0,\dots,\sigma_i\}$ 
before touching the ground is typically
$\sigma_i^2$ for a SRW, but only if
the walk has good chances to cross the whole slab.
Our key attachement estimate follows.
\bl{lem-DD-low}
Consider diffusive deposition. Let $\sigma$
be a configuration such that $|\sigma|<N/2$. Then
there exists a positive constant $\kappa_D$ such that 
\be{DD-lower-ER}
P\big(\text{ Explorer attaches pile } i\ \big|\ \sigma\big)
\ge \kappa_D \frac{\sigma^2_i+1}{N}\times 
\exp\Big(-\frac{3}{N} \big(\sum_{j=1}^N \sigma_j(\sigma_j+1)\big)
\big(1+O(\frac{|\sigma|}{N}\big)\Big).
\ee
\el
In a sense \reff{DD-lower-ER} and \reff{DD-partition}
are saying opposite things: the former
inequality tells how easy it is to get trapped, whereas
the latter tells how easy it is to reach the ground.

Imagine a regime where $\sum_i \sigma_i^2\gg N$. In view of 
\reff{DD-partition} the probability of hitting the ground would be 
small, and very likely the walk would not go below $H$, 
where $H$ is such that
\be{def-H}
\sum_i\big(\sigma_i-H\big)_+^2\sim N.
\ee
In other words, $H$ of \reff{def-H} would play the role of an effective
ground.
We then replace \reff{DD-lower-ER} by the following estimate.
\bc{cor-effective} In the diffusive case, and for any positive $H$,
\be{threshold-DD}
P\big(\text{ Explorer attaches pile }i\big|\sigma\big)
\ge \kappa_D \frac{(\sigma_i-H)^2_+}{N}\times
\exp\big(-\frac{3}{N} \sum_{x=1}^N \big(\sigma_x-H\big)_+^2
\big(1+O(\frac{|\sigma|}{N}\big)\big).
\ee
\ec

\paragraph{Sketch.} We wish to sketch heuristically the reason
why $\log(N)$ is the critical height at which a monopole forms.
We fix a given column, say column 1, and we estimate the number of
explorers needed to produce a given height. Lemma~\ref{lem-DD-low}
allows us to bound this number by a sum of independent geometric
variables, for which we know everything. Indeed, introduce $\tau_1$ the
number of explorers needed so that the height of our distinguished 
site reaches height 1, that is $\tau_1:=\inf\{n>0:\ \sigma_1(n)=1\}$. 
By induction, for any integer $h$ knowing $\tau_{h}$ we define
$\tau_{h+1}:=\inf\{n>0:\ 
\sigma_1(\tau_{h}+n)-\sigma_1(\tau_{h})=1\}$. Assume now
that we are in a regime where hitting the ground is likely. The
estimate \reff{DD-lower-ER} means that for any integers $h,n$
\be{main-geo}
P\big(\tau_{h+1}>n|\ \tau_1,\dots,\tau_h\big)\le
\big(1-\kappa_D\frac{h^2}{N}\big)^n.
\ee
Let us now introduce independent geometric variables $\{\tilde \tau_h,
\ h\ge 1\}$ with $E[\tilde \tau_{h+1}]=N/(\kappa_D h^2)$. Then,
we will show that for any height $H$
\be{geo-comp1}
P\big(\sum_{i=1}^H \tau_i\le X\big)
\ge P\big(\sum_{i=1}^H \tilde \tau_i\le X\big).
\ee
Recall that $\{\sum_{i=1}^H \tau_i\le X\}$ means that $X$ explorers
produce a column of height $H$ at site 1. Now, we want to find $X$ such that
a given height $H$ is likely to be reached. This would be the case
if the probability that any distinguished site reaches height $H$ 
is above $1/N$. Thus, we look for $X$ such that
\be{geo-2}
P\big(\sum_{i=1}^H \tilde \tau_i\le X\big)\sim \frac{1}{N}.
\ee
Let us write $X$ as
$N/f(N)$, and try to guess the size of $f(N)$ which produces a 
monopole. Note also that for any $H>f(N)$, 
\[
E\big[ \sum_{i=f(N)}^H \tilde \tau_i\big]
=\sum_{i=f(N)}^H \frac{N}{\kappa_D i^2}\le
\frac{N}{\kappa_D f(N)}.
\]
Thus, $\{\sum_{i=1}^H \tilde \tau_i\le N/f(N)\}$
imposes a constraint only on the first $f(N)$ variables
in \reff{geo-2}. We then have to estimate $f(N)$ such that
\[
P\big(\sum_{i=1}^{f(N)} \frac{\tilde \tau_i}
{N}\le \frac{1}{f(N)}\big)\sim 
\frac{1}{N}.
\]
We use now that $\{\tilde \tau_i\}$ are independent geometric
variables
\be{geo-4}
\begin{split}
P\big(\tilde\tau_1+\dots+\tilde\tau_{f(N)}\le \frac{N}{f(N)}\big)\ge&
\prod_{i=1}^{f(N)}P\big(\tilde \tau_i\le \frac{N}{f^2(N)}\big)\ge
\prod_{i=1}^{f(N)}\Big(1-\big(1-\frac{\kappa_D i^2}{N}\big)^{N/f^2(N)}
\Big)\\
\ge & \prod_{i=1}^{f(N)}\Big(1-\exp(-\kappa_D \frac{i^2}{f^2(N)})\Big)
\sim
\prod_{i=1}^{f(N)}\Big(\kappa_D \frac{i^2}{f^2(N)}\Big)\\
= & \Big(\kappa_D \frac{1}{f^2(N)}\Big)^{f(N)}(f(N)!)^2.
\end{split}
\ee
Now, using Stirling's formula, we obtain
\be{geo-5}
P\big(\tilde\tau_1+\dots+\tilde\tau_H\le X\big)\ge
\Big( \kappa_D \frac{1}{e^2}\Big)^{f(N)}=
\exp\Big(-\log\big(\frac{e^2}{\kappa_D}\big)\times f(N)\Big).
\ee
Thus, if $f(N)$ is of order $\log(N)$, it is likely that 
one monopole forms.
\section{Comparison with Urns}
\label{sec-urns}
In this section, we establish 
a coupling between our deposition processes and simpler ones
which preserves a natural order on ordered configurations,
to be defined below.

We consider growth evolution on $\N^N$ such that at
each unit time we add a unit height to a configuration,
say $\eta$, at a given site, say $i$, with a probability $p_i(\eta)$ 
which depends only on the value $\eta_i$, and on
the unordered set $\{\eta_j,j\not= i\}$.
In this case, it is useful to reorder the indices through a permutation
of the indices to obtain configurations whose heights are in
decreasing order.
We call $p=\{(p_1(\eta),\dots,p_N(\eta)),\ \eta\in \N^N\}$ the law
of the growth process. 

\subsection{Comparing evolutions}
\label{s:com-ev}
It will be important to compare configurations with
the same number of explorers.
Our main results are the following. 

\bp{lem-coupling} 
Consider two processes 
$t\mapsto\eta(t)$ and $t\mapsto\sigma(t)$ on $\N^N$ evolving, 
respectively, according to the laws 
$p=\big(p_1(\cdot),\dots,p_N(\cdot)\big)$ and
$q=\big(q_1(\cdot),\dots,q_N(\cdot)\big)$.
Assume that for any $\eta,\sigma\in\O_N$ such that 
$\eta\prec\sigma$ we have that 
\begin{equation}
\label{lemco}
\forall k=1,\dots,N\qquad \sum_{i=1}^k p_i(\eta)\le
\sum_{i=1}^k q_i(\sigma).
\end{equation}
Then, the process $\sigma(t)$ is more monopolistic 
than $\eta(t)$,i.e. there is a coupling between the two processes 
such that $\sigma(t)\succ\eta(t)$ for any $t$. 
\ep

\bl{lem-comp}
Let $\{p_1,\dots,p_N\}$ and 
$\{q_1,\dots,q_N\}$ two sets of positive
numbers both summing up to 1. Assume that
\be{hyp-1}
\frac{q_1}{p_1}\ge \frac{q_2}{p_2}\ge\dots\ge \frac{q_N}{p_N}.
\ee
Then, for any $k=1,\dots,N$, we have that
\be{ineq-main}
\sum_{i=1}^k q_i\ge \sum_{i=1}^k p_i.
\ee
\el
\bpr
The proof is by induction on $N$. Assume the Lemma is true
with $N-1$ sets of positive numbers, and define the renormalized
$N-1$ numbers
\be{def-renormalize}
\tilde p_i=p_i\times \frac{1}{1-p_N},\quad\text{and}\quad
\tilde q_i=q_i\times \frac{1}{1-q_N}.
\ee
The induction hypothesis states that for $k=1$ to $N-1$
\[
\sum_{i=1}^k \tilde q_i\ge \sum_{i=1}^k \tilde p_i.
\]
In other words,
\be{step-4}
\sum_{i=1}^k q_i\ge \frac{(1-q_N)}{(1-p_N)}\sum_{i=1}^k p_i.
\ee
The question is whether $1-q_N\ge 1-p_N$ or $q_N/p_N\le 1$
which follows from \reff{hyp-1} since
\[
1=\sum_{i=1}^N p_i\frac{q_i}{p_i}\ge \sum_{i=1}^N p_i\frac{q_N}{p_N}=
\frac{q_N}{p_N}.
\]
\epr

As a corollary of Proposition~\ref{lem-coupling} and 
Lemma~\ref{lem-comp} we have the following result.
\bc{cor-coupling} 
With the notation of Proposition~\ref{lem-coupling},
assume that for any $\eta,\sigma\in\O_N$ such that 
$\eta\prec\sigma$ we have that 
\[
\frac{p_i(\eta)}{p_{i+1}(\eta)}\le 
\frac{q_i(\sigma)}{q_{i+1}(\sigma)}
\;\textrm{ for }i=1,\dots,N-1.
\]
Then, there is an order preserving coupling between $\eta(t)$
and $\sigma(t)$.
\ec

In order to prove Proposition~\ref{lem-coupling} we need 
some notation and some simple observations.
We define the action $\A_j:\N^N\to\N^N$ 
of adding one explorer to site $j$:
$(\A_j\eta)_i=\eta_i+\delta_{i,j}$.
Note that $\A_j$ does not leave $\O_N$ invariant.

Assume that $\eta\in \O_N$ and define
\be{def-increase}
I(\eta)=\acc{i\in \{2,\dots,N\}:\ \eta_{i-1}>\eta_i}\cup\{1\}.
\ee
Also, for $i\in \{1,\dots,N\}$, let 
$d(\eta,i)=\max I(\eta)\cap\{1,\dots,i\}$.
In other words, $d(\eta,i)$
is the last position of a height decrease up to position $i$. 
Note that for $\eta\in \O_N$
we have 
\be{action-T}
\overline{\A_i\eta}
=\A_{d(\eta,i)}\eta\in \O_N.
\ee
For $\eta\in \O_N$, note that 
if $i\le j$, then $d(\eta,i)\le d(\eta,j)$. Also, 
if $i\le j$, then $\A_j\eta\prec\A_i\eta$.

The main observation about ordering is the following.
\bl{lem-order}
Assume that $\eta,\sigma\in\O_N$ with $\eta\prec \sigma$.
If $i\le j$, then $\overline{\A_j\eta}\prec \overline{\A_i\sigma}$.
\el
This lemma is based on the following simple observation.
\bl{lem-step}
Assume that $\eta,\sigma\in \O_N$ with $\eta\prec \sigma$.
If some integers $i<j$ satisfying $d(\eta,j)\le i$, 
are such that
\be{hyp-equal}
\sum_{k=1}^i \eta_k=\sum_{k=1}^i \sigma_k,
\ee
then $i\ge d(\sigma,j)$.
\el

\par\noindent
\textit{Proof of Lemma~\ref{lem-step}.\/}
Assume for a moment that $i>1$. Since $\eta\prec\sigma$, we have
\be{order-1}
\sum_{k=1}^{i-1} \eta_k\le \sum_{k=1}^{i-1} \sigma_k,
\quad\text{and}\quad
\sum_{k=1}^{i+1} \eta_k\le \sum_{k=1}^{i+1} \sigma_k.
\ee
Rewrite now \reff{hyp-equal} as
\[
\sum_{k=1}^{i-1} \eta_k+\ \eta_i=\sum_{k=1}^{i-1} \sigma_k+\ \sigma_i,
\quad\text{and}\quad
\sum_{k=1}^{i+1} \eta_k-\ \eta_{i+1}=\sum_{k=1}^{i+1} \sigma_k-\ 
\sigma_{i+1}.
\]
Using \reff{order-1}, we have both that $\sigma_i\le \eta_i$, and
$\eta_{i+1}\le \sigma_{i+1}$. Since $\sigma\in \O_N$, 
we have $\sigma_{i+1}\le \sigma_i$.
Now, if $i=1$, we have $\sigma_i=\eta_i$ so that $\sigma_i\le \eta_i$
is again true. So that we reach
\be{order-2}
\eta_{i+1}\le \sigma_{i+1}\le \sigma_i\le \eta_i.
\ee
Now $d(\eta,j)\le i<j$ means that $\eta_i=\eta_{i+1}=\dots=\eta_j$,
and with \reff{order-2} this implies that $\sigma_i=\sigma_{i+1}$, and
by induction, we reach that
\be{order-3}
\eta_i=\sigma_i=\sigma_{i+1}=\dots=\sigma_j.
\ee
These last equalities mean that $i\ge d(\sigma,j)$.
\qed

\par\noindent
\textit{Proof of Lemma~\ref{lem-order}.\/}
To simplify notation assume $\eta,\sigma\in \O_N$.
We have already observed that for $i\le j$, $\A_j\eta\prec \A_i \eta$. Thus,
we only need to prove that $\A_j\eta\prec \A_j\sigma$.
If $d(\eta,j)=j$, then $d(\eta,j)\ge d(\sigma,j)$, and the result is 
obvious. Assume henceforth that $d(\eta,j)<j$. If for all $k=d(\eta,j),
\dots,j-1$, we have that
\[
\sum_{i=1}^{k} \eta_i<\sum_{i=1}^{k} \sigma_i,
\]
then the result is also obvious.
In the opposite case, let $k$ in $[d(\eta,j),j[$, be the
first index for which we have
\[
\sum_{i=1}^{k} \eta_i=\sum_{i=1}^{k} \sigma_i,
\]
then, Lemma~\ref{lem-step} implies that $k\ge d(\sigma,j)$,
and the lemma follows.
\qed

\par\noindent
\textit{Proof of Proposition~\ref{lem-coupling}.\/}
By way of induction,
assume that up to time $t$, we have $\eta(t)\prec \sigma(t)$.
Draw a uniform random variable $U$ in
$[0,1[$, and define two random variables $J,J^*$ as follows:
\begin{itemize}
\item[--] if $U\in [p_1(\eta(t))+\dots+p_{i-1}(\eta(t)),p_1(\eta(t))
+\dots+p_{i}(\eta(t))[$ then $J=i$ (we set $p_0=0$);
\item[--] if $U\in [q_1(\sigma(t))+\dots+q_{j-1}(\sigma(t)),q_1(\sigma(t))
+\dots+q_{j}(\sigma(t))[$ then $J^*=j$.
\end{itemize}
Then, \reff{lemco} implies that $J\ge J^*$. We set
\be{coupling-1}
\eta(t+1)=\overline{\A_J\eta(t)}\qquad\text{and}\quad
\sigma(t+1)=\overline{\A_{J^*}\sigma(t)}.
\ee
Then, Lemma~\ref{lem-order} yields that $\eta(t+1)\prec \sigma(t+1)$.
\qed 

\subsection{Comparing Polya's Urn with Random Allocation}
\label{s:com-pr}
By {\it random allocation}, we mean repeated draws of
one out of $N$ colors, labelled from 1 to N, uniformly
at random. In other words, at each draw, the probability
to pick up color $i$ is $1/N$. 
The law for Polya's urn and random allocation 
are denoted respectively $q^P$ and $q^U$ with 
\begin{displaymath}
\forall \sigma\in\N^N,\qquad
q_i^P(\sigma)=
\frac{\sigma_i+1}{N+\sum_{i\le N} \sigma_i}
\qquad\textrm{ and }\qquad
q_i^U(\sigma)=\frac{1}{N}
\;\;.
\end{displaymath}

\bl{lem-polya-uni}
Polya's urn with $N$ colors is more monopolistic than
random allocation of $N$ colors.
\el

\par\noindent
\textit{Proof.\/}
Note that for any $\sigma,\eta\in\O_N$ such that $\eta\prec\sigma$
\be{main-free}
\frac{q_i^P(\sigma)}{q_{i+1}^P(\sigma)}\ge
\frac{q_i^U(\eta)}{q_{i+1}^U(\eta)}
\Longleftrightarrow
\frac{\sigma_{i}+1}{\sigma_{i+1}+1}\ge 1,
\ee
which clearly holds since $\sigma\in\O_N$. 
Thus, Corollary~\ref{cor-coupling} implies the lemma.
\qed

\par\noindent
\textit{Proof of Proposition~\ref{theo-polya}.\/}
Note that if $\eta\succ\sigma$ and $|\eta|=|\sigma|$, we have 
\begin{displaymath}
\forall k=1,\dots,N,\qquad \sum_{i=1}^k q_i^P(\eta) =
\frac{k+\sum_{i=1}^k \eta_i} {N+\sum_{i=1}^N \eta_i}
\ge \frac{k+\sum_{i=1}^k \sigma_i} {N+\sum_{i=1}^N \sigma_i}
= \sum_{i=1}^k q_i^P(\sigma).
\end{displaymath}
This establishes Proposition~\ref{theo-polya}
saying that Polya's evolution with $N$ colors preserves the order.
\qed

\subsection{Comparing deposition models with Polya's urn}
\label{s:com-dp}
Recall the definition of ballistic and diffusive 
deposition given in Section~\ref{sec-intro}.
Denote their law, respectively, by $p^B$ and $p^D$.
We show that both ballistic and diffusive deposition
are more monopolistic than Polya's urn, which is 
one of the statements of Theorem~\ref{theo-comparison}.

Assume for a moment the following lemma.
\bl{lem-BDP}
For any $k=1,\dots,N$ and $\eta\in\O_N$
\be{BDP}
\sum_{i=1}^k p_i^{B}(\eta) \ge \sum_{i=1}^k q_i^P(\eta) 
\;\;\textrm{ and }\;\;
\sum_{i=1}^k p_i^{D}(\eta) \ge \sum_{i=1}^k q_i^P(\eta) \;\;.
\ee
\el
By Proposition~\ref{theo-polya}, Lemma~\ref{lem-BDP}, and 
Proposition~\ref{lem-coupling} we have that both
ballistic and diffusive deposition are more monopolistic than 
Polya's urn. 

To state a preliminary simple observation, we need more
notation. For $\eta\in\O_N$, let
$p^{B,D}_{i,k}(\eta)$ be the probability that the explorer hits
site $i$ at height $k$, for $k\in \N$.
\begin{lemma}
\label{t:bd000}
For any $i,\in \{1,\dots,N\}$ and $\eta(i)\ge k>k'\ge 0$ we have
\be{order-slice}
p^{B}_{i,k'}(\eta) < p^{B}_{i,k}(\eta),\quad\text{and if
$\eta(j)\ge k$}
\quad p^{B}_{i,k}(\eta)=p^{B}_{j,k}(\eta)
\ee
Similarly, 
\be{order-slice02}
p^{D}_{i,k'}(\eta) < p^{D}_{i,k}(\eta),\quad\text{and if
$\eta(j)\ge k$}
\quad p^{D}_{i,k}(\eta)=p^{D}_{j,k}(\eta)
\ee
\end{lemma}

\par\noindent
\textit{Proof of Lemma~\ref{lem-BDP}.\/}
We consider first the ballistic case.
In view of Lemma~\ref{lem-comp}, we need to show that
\be{step-5}
\frac{p_i^{B}(\eta)}{p^{B}_{i+1}(\eta)}
\ge
\frac{\eta_i+1}{\eta_{i+1}+1}
=
\frac{q_i^P(\eta)}{q^P_{i+1}(\eta)}
\;\;.
\ee
In order to prove \eqref{step-5}, we need to show that
\be{step-7}
\frac{\sum_{k=1}^{\eta_i}p^{B}_{i,k}(\eta)}{\eta_i+1}
\ge
\frac{\sum_{k=1}^{\eta_{i+1}}p^{B}_{i+1,k}(\eta)}{\eta_{i+1}+1}
\ee

By Lemma~\ref{t:bd000}, this inequality has the structure
\begin{displaymath}
\frac{a_1+\dots+a_n}{n}\ge\frac{a_{m+1}+\dots+a_n}{n-m}
\end{displaymath}
for $n>m\ge1$ and $a_1\ge a_2\ge\cdots\ge a_n$. The validity of such 
an inequality is immediate once we let $\mu=(a_{m+1}+\dots+a_n)/(n-m)$, 
note $a_1\ge\cdots\ge a_{m}\ge\mu$, 
and write
\begin{displaymath}
\frac{a_1+\dots+a_n}{n}
=
\frac{a_1+\dots+a_{m}+\mu(n-m)}{n}
\ge
\frac{\mu m+\mu(n-m)}{n}
=\mu
\end{displaymath}

Finally, by using \eqref{step-5} and Lemma~\ref{lem-comp} the 
first of equations \eqref{BDP} follows immediately. 
The diffusive case can be treated in the same way. 
This completes the proof 
of the lemma. 
\qed

\par\noindent
\textit{Proof of Lemma~\ref{t:bd000}.\/}
First we prove the lemma for the ballistic case. 
The lemma follows, since
\begin{displaymath}
p^B_{i,k}(\eta)
=
\frac{1}{N}
P(\textrm{the explorer survives till }k+1)
\end{displaymath}
and 
$P(\textrm{the explorer survives till }k)$ is an increasing function 
of $k$. 

Now, we consider diffusive deposition. 
First, assume $k-k'$ an even number.
To each path $s'=\{(x'_j,z'_j)\}_{j=1,...,n}$ hitting $\eta$ in $i$ at height $k'$ 
in a time $n$ we associate uniquely 
a path $s=\{(x_j,z_j)\}_{j=1,...,n}$, of the same length 
and, therefore, the same probability, hitting $\eta$ in $i$ at height $k$.
Thus, the lemma follows since there exist other paths ending in $i$ 
at height $k$. 
Given the path $s'$, we construct $s$ in the following way.
Call $n_1$ the time of last passage of $s'$ through the intermediate height 
$H=(k+k')/2$. $s$ is equal to $s'$ up to time $n_1$ 
while after $n_1$ it uses opposite height increments than the original $s'$,
i.e. $z_{j+1}-z_j=-(z'_{j+1}-z'_j)$ for all $n_1\le j< n$
keeping the same horizontal increments. We therefore obtain a
path ending in $i$ and height $k$. Note that such a path
avoids $\eta$ until it hits site $i$ at height 
$k$, because
$\eta$ is a union of columns.
If $k-k'$ is an odd number we do a similar construction
but we have to associate a set of paths $s'$ of lenght $n$ to a single path
$s$ of lenght $n-1$. The set is obtained considering together
all the paths $s'$ coinciding everywhere but the component $x_{n_1+1}$,
where $n_1$ is now the last hitting time of the level $\frac{k+k'+1}{2}$ of the vertical process.
Then the path $s$ hitting $\eta$ on site $i$ at height $k$ in a time $n-1$
coicide with the paths $s'$ up to time $n_1$, say $x_j=x'_j,\; z_j=z'_j$ for any $j\le n_1$
and is specular to them after $n_1+1$, that is $x_j=x'_{j+1},\; z_j-z_{j-1}=-(z'_{j+1}-z'_j)$ for any $j=n_1+1,..., n-1$ so that
$x_{n-1}=i$ and this is the first hitting to $\eta$.
Clearly the probability of $s$ is larger
than or equal to the sum of the probabilities
of the paths $s'$, since $s$ is one step shorter, 
and the sum on $x_{n_1+1}$ is 
done only on $N-\zeta_{n_1+1}(\eta)$ sites.

%
\qed
\section{Estimating Unit Growth}
\label{sec-unit}
In this section we discuss how heights grow.
We consider the random walk $S_n=(X_n,Z_n)$, and for an integer
$k$, we call $H_k$ the first time the walk reaches height $k$. 
In other words,
\be{def-Hk}
H_k=\inf\{n\ge 0:\ Z_n=k\}.
\ee
Since the $X$-component is uniform on the base, giving the configuration
$\sigma$, the ordered one $\bar\sigma$, 
or the height occupation $\zeta$ (defined in \eqref{heightocc})
is equivalent, and we use $P_g(\sigma)$ or $P_g(\zeta)$ indifferently
to denote the probability an explorer hits the ground.
Lemma~\ref{lem-DD} is obtained as a simple application of Jensen's
inequality whereas Lemma~\ref{lem-DD-low} requires
Kesten-Kozlov-Spitzer representation of the local times
\cite{KKS}.

\paragraph{Proof of Lemma~\ref{lem-DD}.}
For an integer $k\le \bar\sigma_1$
let $l(k)$ be the number of visits of height $k$ by the
random walk before $H_0$. We have the represenation
\be{represent-Z}
P_g(\zeta)=E\big[ \prod_{k=1}^{\bar\sigma_1}\big(1-\frac{\zeta_k}{N}\big)^
{l(k)}\big]=E\big[ \exp\Big(\sum_{k=1}^{\bar\sigma_1} l(k)
\log\big(1-\frac{\zeta_k}{N}\big)\Big).
\ee
Our hypothesis $|\sigma|<N/2$ implies that for $k\ge 1$, 
$\zeta_k\le \zeta_1\le N/2$,
and since $\log(1-x)\ge -x-x^2$ for $0\le x\le 1/2$, we have using
Jensen's inequality
\be{main-Z-lb}
P_g(\zeta)\ge \exp\Big(-\sum_{k=1}^{\bar\sigma_1}E[l(k)]
\big(\frac{\zeta_k}{N}+\frac{\zeta_k^2}{N^2}\big)\Big).
\ee
Now note that $E[l(k)]=2k$. Indeed, the height of the random walk being a 
simple random walk on $\N$, we have for $k\le \bar\sigma_1$,
by conditionning on the first step
\be{ruin-1}
\begin{split}
E[l(k)|Z_0=\bar\sigma_1]=&E[l(k)|Z_0=k]=1+\frac{1}{2}
\big(E[l(k)|Z_0=k+1]+E[l(k)|Z_0=k-1]\big)\\
=& 1+\frac{1}{2}
\big(E[l(k)|Z_0=k]+P(H_k<H_0|Z_0=k-1)E[l(k)|Z_0=k]\big)\\
=&1+\frac{1}{2} \big(E[l(k)|Z_0=k]+(1-\frac{1}{k})E[l(k)|Z_0=k]\big).
\end{split}
\ee
The equality $E[l(k)]=2k$ for $k\le \bar\sigma_1$ follows at once.
Note now that
\[
\sum_{k\ge 1} 2k\zeta_k=\sum_{i=1}^N \sigma_i(\sigma_i+1).
\]
Finally, \reff{DD-partition} follows as we note that
\[
\sum_{k\ge 1} 2k
\frac{\zeta_k^2}{N^2}\le \frac{\zeta_1}{N}\sum_{k\ge 1} 
2k\frac{\zeta_k}{N}
\le 2\frac{|\sigma|}{N}\sum_{i=1}^N \frac{\sigma_i(\sigma_i+1)}{N}.
\]
\qed
\paragraph{On Kesten-Kozlov-Spitzer representation.}
Let $u(h)$ be the number of {\it up-crossings} of height $h$
before touching the base. In other words, we define $u(0)=0$ and for
$h>0$
\be{def-u}
u(h)=\sum_{i=1}^{H_0} \ind_{(Z_{i-1},Z_i)=(h,h+1)}.
\ee
Similary, {\it down-crossings} of height $h$ correspond to jumps
from $h$ to $h-1$ before time $H_0$. 
One way to realize the random walk $n\mapsto Z_n$
is to assign the sequence of up and down-crossings on each height.
Thus, we consider $\{\{\xi_i^k,\ i\in \N\},\ k\in \N\}$
a collection of i.i.d. geometric variables,
with law $P(\xi=n)=1/2^{n+1}$ for $n\in \N$. 
Now, the sequence of up and down-crossings 
at height $k$ is as follows: $\xi_0^k$ up-crossings, 
then one down-crossing, then
$\xi_1^k$ up-crossings, the one down-crossing, 
then $\xi_2^k$ up-crossings...
and so on and so forth. The key observation is that each $\xi_i^k$,
for $i\ge 1$, is preceded by an up-crossing of the height $k-1$.
In other words, 
\be{rec-u}
u(k)=\xi_0^k+\sum_{j=1}^{u(k-1)}\xi_j^k,\quad\text{and}\quad
u(1)=\xi_0^1.
\ee

We set $\kG(h)= \sigma( \xi_i^k,\ k\le h,\ i\in \N)$ the $\sigma$-field 
representing the choices of moves on the first $h$ heights.
Kesten-Kozlov-Spizter representation expresses the local times
of $Z$ in terms of the $u$. Thus, if $l(k)$ represents the number
of visits of height $k$ before $H_0$, for a walk with starting level
above $\bar\sigma_1$, then 
\be{def-KKS}
\forall k\ge 1,\quad l(k)=u(k)+u(k-1)+1.
\ee
Then, with notation $x_i=1-\zeta_i/N$
\be{Z-1}
\begin{split}
P_g(\zeta)=& E\big[\prod_{k=1}^{\bar\sigma_1}x_k^{l(k)}\big]=
E\big[\prod_{k=1}^{\bar\sigma_1-1}x_k^{l(k)}
x_{\bar\sigma_1}^{u(\bar\sigma_1)+u(\bar\sigma_1-1)+1}\big]\\
=&E\big[ \prod_{k=1}^{\bar\sigma_1-1} x_k^{l(k)}x_{\bar\sigma_1}^{u(\bar\sigma_1-1)+1}
E\big[x_{\bar\sigma_1}^{\xi_0^{\bar\sigma_1}+\dots+\xi_{u(\bar\sigma_1-1)}^{\bar\sigma_1}}
\big| \kG(\bar\sigma_1-1)\big]\big].
\end{split}
\ee
Since $E[z^\xi]=1/(2-z)$, we have
\be{Z-2}
E\big[x_{\bar\sigma_1}^{\xi_0^{\bar\sigma_1}+\dots+\xi_{u(\bar\sigma_1-1)}^{\bar\sigma_1}}
\big| \kG(\bar\sigma_1-1)\big]=\frac{1}{(2-x_{\bar\sigma_1})^{1+u(\bar\sigma_1-1)}}.
\ee
We set $a(k)=0$ for $k\ge\bar\sigma_1$, whereas for any $k< \bar\sigma_1$
\be{Z-3}
e^{-a(k-1)}=\frac{x_k}{2-x_k e^{-a(k)}},
\ee
so that
\be{Z-4}
P_g(\zeta)=
e^{-a(\bar\sigma_1-1)} E\big[ \prod_{k=1}^{\bar\sigma_1-1}x_k^{l(k)}
e^{-a(\bar\sigma_1-1)u(\bar\sigma_1-1)}\big],
\ee
and by induction, we obtain
\be{Z-5}
P_g(\zeta)=\exp\big(-\sum_{k\ge 0} a(k)\big).
\ee
Note that \reff{Z-3} reads for $1\le k\le \bar\sigma_1$
\be{Z-6}
e^{-a(k)}+e^{a(k-1)}=\frac{2}{x_k},\quad\text{and}\quad
e^{a(k-1)}-e^{a(k)}=e^{-a(k+1)}-e^{-a(k)}+\frac{2}{x_k}-\frac{2}{x_{k+1}},
\ee
and $a(k-1)\ge a(k)\ge 0$ follows by induction from \reff{Z-3}. 
%
%
%
Inequality \reff{DD-partition} implies that
\be{Z-9}
\sum_{k\ge 0} a(k) \le \frac{2}{N}\sum_{j\ge 1} j\zeta_j+
\frac{2}{N^2}\sum_{j\ge 1} j\zeta^2_j.
\ee
\paragraph{Proof of Lemma~\ref{lem-DD-low}.\/}
An explorer settling on the pile at site $i$, hits the $i$-th pile
at a height between 1 and $\sigma_i$. Knowing that it settles 
at height $k$, it has chance $1/\zeta_k$ to settle on $(i,k)$
since we are on the complete graph. We underestimate the 
probability of settling on $\sigma_i$, if we only consider
trajectories hitting only one of the $\{\zeta_1,\zeta_2,\dots\}$
before $H_0$. Thus,
\be{p1-1}
\begin{split}
P\big(\text{Explorer attaches pile $i$}\big|\ \sigma\big)
\ge&\sum_{h= 1}^{\sigma_i}\frac{1}{\zeta_h}
E\big[\big(1-x_h^{l(h)}\big)\prod_{j\not= h} x_j^{l(j)}\big]\\
\ge & \sum_{h=1}^{\sigma_i}\frac{1}{\zeta_h}
\Big(E\big[\prod_{j\not= h} x_j^{l(j)}-P_g(\zeta)\Big).
\end{split}
\ee
Fix $h>0$, and write $\zeta^h$ for the height occupation such that
\[
\forall k\not= h,\quad\zeta^h_k=\zeta_k,\quad\text{and}\quad
\zeta^h_h=0.
\]
We rewrite \reff{p1-1} in terms of the function $P_g$ as follows
\be{p1-2}
P\big(\text{Explorer attaches pile $i$}\big|\ \sigma\big)
\ge P_g(\zeta)\sum_{h= 1}^{\sigma_i}\frac{1}{\zeta_h}
\big( \frac{P_g(\zeta^h)}{P_g(\zeta)}-1\big).
\ee
For a given $h\le \bar\sigma_1$, we now study the ratio 
$P_g(\zeta^h)/P_g(\zeta)$.
As in \reff{Z-5} in the previous paragraph, we write
$P_g(\zeta^h)=\exp(-\sum_{k\ge 0} \tilde a(k))$, with $\tilde a$
satisfying the relation \reff{Z-3} with $\zeta^h$ in place of $\zeta$.
In other words,
\be{def-tildea}
\forall k\ge h,\quad \tilde a(k)=a(k),\quad 
\text{and}\quad\exp(-\tilde a(h-1))= \frac{1}{2-e^{-a(h)}}.
\ee
For $k\le h-1$,
\be{p1-cond}
e^{\tilde a(k-1)}+e^{-\tilde a(k)}=\frac{2}{x_k}=
e^{a(k-1)}+e^{-a(k)}.
\ee
We set $\delta_k=a(k)-\tilde a(k)$, and from \reff{p1-cond}, we have
that $\delta_k\ge 0$. In terms of $\delta_k$, \reff{p1-cond} reads
for $k<h-1$
\be{p1-3}
\exp(\delta_{k-1})-1=e^{-\tilde a(k-1)-a(k)}\big(\exp(\delta_{k})-1\big),
\ee
whereas for $k=h-1$ we have
\be{p1-bis}
\exp(\delta_{h-1})-1=\frac{2(1-x_h)}{x_h(2-\exp(-a(h))}\le 2.
\ee
Equality \reff{p1-bis} implies (since $a(.)\le 2$ by \reff{p1-cond})
that for some constant $\kappa_D$
\be{p1-last}
\delta_{h-1}\ge \kappa_D(1-x_h)=\kappa_D\frac{\zeta_h}{N}.
\ee
We deduce that $\delta_{k-1}\le \delta_k$ and since $x\mapsto (e^x-1)/x$
is increasing, \reff{p1-3} implies
\be{p1-4}
\frac{\delta_{k-1}}{\delta_k} \ge \exp\big(-\tilde a(k-1)-a(k)\big)\ge
\exp\big(-a(k-1)-a(k)\big)\ge \exp\big(-2a(k-1)\big).
\ee
By induction on \reff{p1-4}, and using \reff{p1-last}, this implies that
for each $h\le \bar\sigma_1$
\be{p1-5}
\begin{split}
\frac{P_g(\zeta^h)}{P_g(\zeta)}-1=&\exp\big(\sum_{k=0}^{h-1}\delta_k\big)-1
\ge \sum_{k=1}^{h-1}\delta_k\\
\ge & \delta_{h-1}\sum_{k=0}^{h-1} \exp\big(-2\sum_{j=k}^{h-1}a(j)\big) \\
\ge & \kappa_D
\zeta_h\frac{(h-1)}{N} \exp\big(-2\sum_{j=0}^{\bar\sigma_1-1}a(j)\big)\\
\ge & \kappa_D\zeta_h\frac{(h-1)}{N} \exp\Big(-4\sum_{j\ge 1}
\big(j\frac{\zeta_j}{N}+j\frac{\zeta^2(j)}{N^2}\big)\Big).
\end{split}
\ee
Since $\zeta_j\le N$, we have $\zeta^2_j/N\le \zeta_j$, and
\reff{p1-5} implies \reff{DD-lower-ER}.
\qed

\br{rem-height} In order to obtain Corollary~\ref{cor-effective},
consider simply the height occupation
$\{\zeta_{k+H},\ k\in \N\}$ and proceed along the
exact same proof. This height occupation corresponds to
$\{(\sigma_i-H)_+,\ i=1,\dots,N\}$.
\er

Given $A>0$, we define the {\it early regime} 
as the following subset of configurations.
\begin{equation}
\label{ear-reggg}
\ER(A)=\{\sigma:\, \sum_{i=1}^N\sigma_i^2\le AN,\; 
\sum_{i=1}^N \sigma_i<\frac{N}{2} \}.
\end{equation}
Define also $\kappa(A)$ to be $\kappa_D\exp(-2A)$.
\bc{cor-ear}
For any $\sigma \in \ER(A)$, and $i\in \{1,\dots,N\}$
\be{hit-ER}
P\big(\text{ Explorer attaches pile }i\big|\sigma\big)
\ge \kappa(A) \frac{\sigma_i^2}{N}.
\ee
\ec
\subsection{Upper Bound}
\label{s:upp}

\bl{lem-DD-upp}
Consider diffusive deposition. Let $i$ be a fixed 
site and $\sigma$ be a configuration such that 
$\sigma_i<\sqrt N$. Then 
\be{DD-upper-ER}
P\big(\text{ Explorer attaches to site } i\ \big|\ \sigma\big)
\le \kappa \frac{(\sigma_i\vee1)^2}{N}\;\;.
\ee
with $\kappa=1+O(N^{-1/2})$.
\el
\br{r}
Comparing lower and upper bound (Lemmas \ref{lem-DD-low} 
and \ref{lem-DD-upp}) on attachment probability we obtain
a good control on this probability for configurations in the early regime,
that is in $\ER(A)$ when $A$ is small.
\er
\smallskip
\par\noindent
\textit{Proof of Lemma~\ref{lem-DD-upp}.\/}
If $\sigma_i=0$ then \reff{DD-upper-ER} is immediate with $\kappa=1$.
If $\sigma_i\ge 1$ the chances an explorer attaches to column $i$, 
in configuration $\sigma$, is smaller than if
all columns distinct from $i$ were set to zero.
This is seen by coupling. First, for configuration $\sigma$,
let $\sigma^i$ denote the configuration where we annihilate
all columns distinct from $i$. In other words
\[
(\sigma^i)_k=0,\quad\text{when }k\not=i,
\qquad\text{and}\quad (\sigma^i)_i=\sigma_i.
\]
Therefore, in $\sigma^i$, the highest column is $i$ with
height $\sigma_i\ge 1$. Now, the event hit column $i$ in $\sigma^i$
is the complement of the event hit the base first. Since 
$P_g(\sigma^i)$ is the probability the explorer hits first the base in 
configuration $\sigma^i$, by Lemma~\ref{lem-DD}, 
we have that for some $\kappa>0$
\begin{displaymath}
\begin{split}
P\big(\text{An explorer attaches to i}\big|\sigma\big)
\le&
P\big(\text{An explorer attaches to i}\big|\sigma^i\big)\\
=& 1-P_g(\sigma^i) \le \frac{\sigma_i^2}{N} (1+O(N^{-1/2}),
\end{split}
\end{displaymath}
which completes the proof.
\qed

\bl{lem-DB-upp}
Consider ballistic deposition. Consider a configuration 
$\sigma\in\N^{G_N}$, then 
\be{DB-upper-ER}
P\big(\text{ Explorer attaches to site } i\ \big|\ \sigma\big)
\le \frac{\sigma_i+1}{N}\;\;.
\ee
\el

\smallskip
\par\noindent
\textit{Proof of Lemma~\ref{lem-DB-upp}.\/}
To get attached to site
$i$, the particle has to survive up to the time it
reaches height $\sigma_i$,
and then at each step--down has a chance $1/N$ to fall on
column $i$ provided it has avoided the other columns.
Thus
\begin{displaymath}
\begin{array}{l}
{\displaystyle
P\big(\text{An explorer attaches to site } i\big|\sigma\big)
}\\
{\displaystyle
\phantom{mm}
=
\bigg[\prod_{h=\sigma_i+1}^{\max \sigma}\Big(1-
\frac{\zeta_h}{N}\Big)\bigg]
  \frac{1}{N}\Big(1+
\Big(1-\frac{\zeta_{\sigma_i}}{N}\Big)+
\Big(1-\frac{\zeta_{\sigma_i}}{N}\Big)
\Big(1-\frac{\zeta_{\sigma_i-1}}{N}\Big)
+\dots\Big)}\\
{\displaystyle
\phantom{mm}
\le  \frac{1}{N}\big(\sigma_i+1\big).}
\end{array}
\end{displaymath}
which completes the proof.
\qed
\subsection{Stochastic domination}
\label{sd}
Both in ballistic and diffusive deposition, we have
a simple upper bound on the probability of attaching
to a given column (see Section~\ref{s:upp}), 
which depends only on its height.
This, in turn, is used to bound the
number of explorers necessary to increase the height
by one unit in terms of a geometric random variable,
for which everything can be computed explicitly.
In other words,
call $\tau_1$ the number of explorers needed so that column
1 reaches height 1. Let $\tau_2$ be the additional number of explorers
needed to reach a height 2, and so on. Note that $\{\tau_1>k\}$ means
that out of $k$ explorers none of them has reached site 1.
These times are used to control the height of column 1 after 
$k$ explorers have been sent, 
\begin{displaymath}
\forall k\ge 1,\quad\forall  H\ge 1 ,\qquad
P_0(\sigma_1(k)>H)= P(\tau_1+\cdots+\tau_H<k).
\end{displaymath}

We need to estimate the sum of the 
$\{\tau_i\}$ with the following general lemma.
\bl{lem-comp-geo}
Let $\tau,T$ be stopping times with respect to a filtration $\{\F_n\}$.
Let $\tau_1:=\tau$ and if $\theta(n)$ is the time-shift by $n$ units,
define inductively
\[
\tau_n:=\tau\circ \theta(\tau_1+\dots+\tau_{n-1}).
\]
Let $\{\tilde \tau_n,\ n\in \N\}$ be independent random variables, which 
are also independent from $\{\tau_n,\ n\in \N\}$.
Assume that for positive integers $\xi\le \xi'$, we have
\be{hyp-order}
P\big(\tau_n>\xi\ ,T>\xi'|\ \F_{\tau_1+\dots+\tau_{n-1}}\big)\ge
P\big(\tilde \tau_n>\xi\big)P(T>\xi').
\ee
Then, for any integer $n$ and $\xi>0$
\be{sum-order}
P\Big(\sum_{i=1}^n \tau_i>\xi\ |\ T>\xi\Big)\ge 
P\Big(\sum_{i=1}^n\tilde \tau_i>\xi\Big).
\ee
Similarly, if instead of \reff{hyp-order} we have
\[
P\big(\tau_n>\xi\ ,T>\xi'|\ \F_{\tau_1+\dots+\tau_{n-1}}\big)\le
P\big(\tilde \tau_n>\xi\big) P\big(T>\xi'\big),
\]
then
\[
P\Big(\sum_{i=1}^n \tau_i>\xi|\ T>\xi \Big)
\le P\Big(\sum_{i=1}^n\tilde \tau_i>\xi\Big).
\]
\el

\smallskip
\par\noindent
\textit{Proof of Lemma~\ref{lem-comp-geo}.\/}
We prove \reff{sum-order} by induction. The step $n=1$ is obvious.
Assume the inequality is true at step $n-1$. 
Using that the variables are positive,
\be{induc-1}
\begin{split}
P\Big(\sum_{i=1}^n& \tau_i>\xi,\ T>\xi\Big)= \sum_{K\le\xi}
P\Big(\sum_{i=1}^{n-1} \tau_i=K,\ T>\xi,\ \tau_n>\xi-K\Big)
+P\Big(\sum_{i=1}^{n-1} \tau_i>\xi,\ T>\xi\Big)\\
= &\sum_{K\le \xi} E\Big[ 
P\Big( \tau_n>\xi-K,\ T>\xi-K\Big|\F_K\Big)
\ind_{\sum_{i=1}^{n-1} \tau_i=K,\ T>K}\Big]
+P\Big(\sum_{i=1}^{n-1} \tau_i>\xi,\ T>\xi\Big)\\
\ge &  \sum_{K\le\xi}E\Big[ P\Big(\tilde \tau_n>\xi-K\Big)
P\big(T>\xi-K\Big| \F_K\big)
\ind_{\sum_{i=1}^{n-1} \tau_i=K,\ T>K}\Big]
+P\Big(\sum_{i=1}^{n-1} \tau_i>\xi,\ T>\xi\Big)\\
= &\sum_{K\le\xi} P\Big(\tilde \tau_n>\xi-K\Big)\times
P\Big(\sum_{i=1}^{n-1} \tau_i=K, T>\xi\Big)
+P\Big(\sum_{i=1}^{n-1} \tau_i>\xi,\ T>\xi\Big)\\
= & P\Big(\sum_{i=1}^{n-1} \tau_i+\tilde \tau_n>\xi,\ T>\xi\Big)
\;\;.
\end{split}
\ee
Now, we can exchange the role played by $\tilde \tau_n$ and
by $\tau_1+\dots+\tau_{n-1}$ in the previous argument,
to use the induction hypothesis. Indeed,
\be{induc-2}
\begin{split}
P\big(\sum_{i=1}^{n-1} \tau_i+\tilde \tau_n>\xi,\ T>\xi\big)=&
\sum_{K=1}^{\xi} P\big(\tilde \tau_n=K\big)\times
P\big(\sum_{i=1}^{n-1}\tau_i>\xi-K,\ T>\xi\big)+
P\big(\tilde \tau_n>\xi\big)\\
\ge& P\big(\sum_{i=1}^n \tilde \tau_i>\xi\big)P(T>\xi).
\end{split}
\ee
The proof of the opposite inequalities follows the same steps.
\qed

\section{Very Early Regime}
\label{sec-VE}
One important step in the cluster growth
is to reach height $\log(N)$. We cover this intermediary step in
the following proposition, even if it 
is included in Theorem~\ref{main-theo}.
\bp{prop-init}
Consider diffusive deposition.
There exist positive constants $b$ and $\gamma$ such that
almost surely, for $N$ large
\[
\max_{i\le N} \sigma_i\big(b\frac{N}{\log N}\big)\ge \gamma\log N.
\]
\ep
This proposition is concerned with what we call the {\it very
early regime}, and it is based on comparison with urn models.

Two scales play an important role in this section: the
time scale $N/\log(N)$, and the space scale $\log(N)/
\log\log(N)$. We therefore introduce notation
\be{def-scale}
T_N=\frac{N}{\log(N)},\quad\text{and}\quad
H_N=\frac{\log(N)}{\log\log(N)}.
\ee
The set of configurations $\s$ with 
maximal height lower than $\gamma\log N$ and
$|\s|\le \beta T_N$ 
is called the {\em very early regime} and is denoted by $\VE$.
In other words,
\[
\VE=\{\sigma: |\sigma|\le 
\beta T_N,\quad \max_x \sigma_x\le \gamma \log(N)\},
\]
and note that 
\[
\forall \s\in \VE,\qquad 
\sum_{x=1}^N \s_x^2\le 
\big(\max_{x\le N} \s_x\big)|\s|\le \beta\gamma N.
\]
If $\tau_A$ is the hitting time of set $A$, we show in this
Section that there are constants $b<\beta$ and $\delta>0$ such that
\be{goal-sec5}
P\big( \tau_{\VE^c}> b T_N\big)\le \exp(-N^\delta).
\ee

\paragraph{Strategy of the proof.} 
We divide time in two periods. In the first, of length $T_N$,
a large number
of columns, of order $N^a$ with $0<a<1$, reach a height $\delta H_N$.
This is the content of Lemma~\ref{lem-vve}, whose
main ingredient is a coupling
between diffusive deposition and random allocation.
In the second period, we use the estimate of
Corollary~\ref{cor-effective} to 
control the growth of these columns together with
sending Poisson waves of explorers to ensure 
the growth of each column independently.

\paragraph{Step 1: Reaching height $H_N$.}
The random allocation evolution is denoted by $n\mapsto \eta(n)$.
Our first lemma deals exclusively with random allocation.
\bl{lem-ra}
For $\alpha\in[\frac{1}{2},1)$,
and $\delta<(1-\a)/2$, we have almost surely, for $N$ large enough 
\begin{equation}\label{ra010}
|\{x: \eta_x(T_N)>\delta H_N\}|
\ge
N^\alpha.
\end{equation}
\el
\bpr
Let $X$ be a Poisson variable of parameter $T_N/2$.
We have
\be{poi-1}
P\big(|\{x: \eta_x(T_N)>\delta H_N\}|<N^\alpha\big)\le
P\big(|\{x: \eta_x(X)>\delta H_N\}|<N^\alpha\big)+
P\big(X>T_N\big).
\ee
Now, $\{\eta_x(X),x=1,\dots,N\}$ are independent Poisson variables
of parameter $1/2\log(N)$. A tedious but simple computation gives
\be{poi-2}
P\big(\eta_1(X)\ge \delta H_N\big)=N^{-2\delta(1+o(1))}.
\ee
Now,
\be{poi-3}
|\{x: \eta_x(X)>\delta H_N\}|=\sum_{x=1}^N \ind_{\{\eta_x(X)\ge
\delta H_N\}},
\ee
and by Bernstein's inequality, for $\alpha<1-2\delta$, and
\be{poi-4}
P\big(|\{x: \eta_x(X)>\delta H_N\}|<N^\alpha\big)\le
\exp(-\frac{N^{1-2\delta}-N^{\alpha}}{2}).
\ee
Note also that from Chebychev's exponential inequality
\be{poi-5}
P\big(X>T_N\big)\le \exp\Big(-\frac{3-e}{2}\,T_N\Big).
\ee
The statement follows.
\epr

In Lemma~\ref{lem-BDP}, we establish that
diffusive deposition, denoted $t\mapsto \s(t)$
is more monopolistic than random
allocation. Thus, there is a coupling such that with
probability 1, when $\s(0)=\eta(0)$, we have for any $t\ge 0$
\[
\forall k\le N,\qquad \sum_{i=1}^k \bar \s_i(t)\ge 
\sum_{i=1}^k \bar \eta_i(t).
\]
Assume now that $\sigma(T_N)\in \VE$, and that
\be{poi-6}
L:=|\{x: \eta_x(T_N)>\delta H_N\}|>L':
=|\{x: \sigma_x(T_N)>\frac{\delta}{2} H_N\}|.
\ee
Then, by our coupling
\be{poi-7}
\sum_{i=1}^L \bar \s_i(T_N)\ge 
\sum_{i=1}^L \bar \eta_i(T_N)\ge \delta H_NL,\quad\text{and}\quad
\sum_{i=L'}^L \bar \s_i(T_N)\le \frac{\delta}{2} H_N L.
\ee
Then
\[
\gamma \log(N) L'\ge \sum_{i=1}^L \bar \s_i(T_N)\ge 
\frac{\delta}{2} H_N L
\Longrightarrow L'\ge \frac{\delta}{2\gamma \log\log(N)} L.
\]
Thus, for any $\alpha>1-2\delta$, we have that $L'>N^{\alpha}$.
We therefore state the result as follows.
\bl{lem-vve}
For $\alpha\in[\frac{1}{2},1)$,
and $\delta<(1-\a)/2$, we have almost surely, for $N$ large enough,
and for the diffusive deposition $t\mapsto \s(t)$,
\begin{equation}
\label{dd-1}
|\{x: \sigma_x(T_N)>\frac{\delta}{2} H_N\}|>N^\alpha.
\end{equation}
\el
By applying Markov's property at time $T_N$
\be{poi-8}
\begin{split}
P\big(\tau_{\VE^c}(\s)&>(b+1) T_N\big)
\le P\big( |\{x:\ \eta_x(T_N)>\delta H_N\}|<N^\alpha\big)\\
&+\sup\Big( P_\s\big(\s(bT_N)\in \VE\big):\ \s\in \VE,\ 
|\{x: \sigma_x>\frac{\delta}{2} H_N\}|>N^\a\Big).
\end{split}
\ee

\paragraph{Step 2: Poisson waves.}
We realize diffusive deposition for times in 
$[T_N,bT_ N]$ by a sequence of
Poisson waves, the $k$-th wave made of $X^{(k)}$
explorers, and $\{X^{(k)},\  k\ge 1\}$ an i.i.d sequence of
Poisson random variables with parameter $x_N$ 
going to infinity with $N$.

Our starting configuration denoted $\s^{(0)}$ satisfies
\[
\s^{(0)}\in \VE,\quad
\Lambda_N:=\{x: \sigma^{(0)}_x>\frac{\delta}{2} H_N\},\quad
\text{and}\quad |\Lambda_N|>N^\a.
\]
Let $\sigma^{(k)}$ be the configuration 
of diffusive deposition starting from $\s^{(0)}$ 
after the $k$-th wave is sent, i.e., 
\begin{displaymath}
\sigma^{(k)}=\sigma\Big(\sum_{\ell=1}^{k}X^{(\ell)}\Big) \qquad
\hbox{ with } \quad |\sigma^{(k)}|=|\s^{(0)}|
+\sum_{\ell=1}^{k}X^{(\ell)}
\end{displaymath}
Define now, using $\kappa_D$ of Corollary~\ref{cor-effective},
\[
\KVE:= \kappa_D\exp(-\gamma \beta ).
\]

We have for the diffusive deposition process
for any $k=1,2,\dots$, if $\sigma^{(k)}\in \VE$ then
\begin{equation}
\label{shoeld00}
p_i(\sigma(t))\ge \frac{K_{\normalfont\textrm{ve}}}{N}
          (\sigma_i^{(k-1)})^2\qquad
\forall t\in [|\sigma^{(k-1)}|,X^{(k)}+|\sigma^{(k-1)}|[,\;
\forall i\in\Lambda_N.
\end{equation}
This immediately follows from Corollary~\ref{cor-effective}.

Consider an auxiliary growth process $\tilde\sigma(t)$ 
which evolves on the sites of $\Lambda_N\cup\{0\}$,
defined iteratively as follows. 
Set $\sigma_0^{(0)}=0$, and for $i\in \Lambda_N$, set
$\tilde\sigma_i^{(0)}=\frac{\delta}{2}H_n$.
Each explorer in the $k$--th wave is attached to site 
$i\in \Lambda_N$ with probability
\begin{equation}\label{pikappa}
p_i^A(k)= \KVE\frac{(\tilde\sigma_i^{(k-1)})^2}{N},
\end{equation}
whereas site 0 grows by one with probability 
$1-\sum_{i\in \Lambda_N} p_i^A(k)$.

The following result is crucial.
\bl{lem-aux}
There exists a coupling between $t\mapsto\sigma(t)$ and 
$t\mapsto\tilde\sigma(t)$ such that if $|\sigma^{(k)}|\le \tau_{\VE^c}$
and $t$ is within the $k$-th wave,
i.e., $ t\in [|\sigma^{(k-1)}|,X^{(k)}+|\sigma^{(k-1)}|[$
\begin{equation}\label{coup}
\forall i\in \Lambda_N,\qquad \sigma_i(t)\ge\tilde\sigma_i(t).
\end{equation}
Moreover, $\{\tilde \sigma_i,\ i\in \Lambda_N\}$ are 
independent, where we used the shorthand notation
$\tilde \sigma_i=\{\tilde \sigma_i^{(k)},k\le \tau_{\VE^c}\}$.
\el
\bpr
The coupling part is simple and we omit it here.
We denote by $\{Y^1_i,\ i\in \Lambda_N\}$ independent 
Poisson variables of parameter $x_Np_{i}^A(1)$. 
We denote by $\kG_1$ the sigma-field generated by $X^{(1)}$ and by
$\{Y^1_i,\ i\in \Lambda_N\}$. We now build $\kG_k$ by induction.
Assume that $\kG_{k-1}$ has been built. Then conditioned on
$\kG_{k-1}$, we fix the height of all sites after the $k-1$-th wave.
Draw a Poisson variable $X^{(k)}$ independent of $\kG_{k-1}$, and
denote by $\{Y^k_i,\ i\in \Lambda_N\}$ the independent 
Poisson variables of parameter $p_{i}^A(k)x_N$ which is
itself $\kG_{k-1}$ measurable. Note that $Y^k_i$ depends only
on $X^{(k)}$ and on the past through 
$Y^{1}_i+\dots+Y^{k-1}_i$, the height of site $i$ after the
$k-1$-th wave. In other words, for any real function $f_i$,
there is a function $\phi_i$ such that
\be{devil-1}
E\big[f_i\big(Y^k_i+\dots+Y^k_i\big)|\kG_{k-1}\big]=
E\big[f_i\big(Y^k_i+\dots+Y^k_i\big)|Y^{1}_i+\dots+Y^{k-1}_i\big]=
\phi_i(Y^{1}_i+\dots+Y^{k-1}_i).
\ee
Note also, that if we integrate only over $X_k$, and for any
real functions $f_i$, for $i\in \Lambda_N$ we have
\be{devil-2}
\begin{split}
E\big[\prod_{i\in \Lambda_N}
f_i\big(Y^k_i+\dots+Y^k_i\big)|\kG_{k-1}\big]=&\prod_{i\in \Lambda_N}
E\big[f_i\big(Y^k_i+\dots+Y^k_i\big)|Y^{1}_i+\dots+Y^{k-1}_i\big]\\
=&
\prod_{i\in \Lambda_N}\phi_i(Y^{1}_i+\dots+Y^{k-1}_i).
\end{split}
\ee
This means that what happens on different sites of $\Lambda_N$
is independent.
\epr

Now, each Poisson wave we send has about $x_N$ explorers,
and we expect to send about $(b-1) T_N/x_N$ waves. 
Recall that for $N$ large, we have a.s. that
$\sigma^{(0)} \in \VE$ and $|\Lambda_N|>N^\alpha$.
Therefore, if $t_N$ denotes the integer part of $bT_N/(2ex_N)$,
and for simplicity $H=\frac{\delta}{2}H_N$
\be{devil-3}
\begin{split}
P_{\sigma^{(0)}}\big(\tau_{\VE^c}>b T_N\big)\le&
P\big(\sum_{k=1}^{t_N}|X^{(k)}|>b 
T_N \big)+P_{\sigma^{(0)}}\big( \max_{i\in \Lambda_N}
\big(H+\sum_{k=1}^{t_N} Y_i^k\big)\ <\ \gamma N\big)\\
\le & e^{-bT_N/2}+
\prod_{i\in\Lambda_N}\Big(1- P_{\sigma^{(0)}}\big(
\big(H+\sum_{k=1}^{t_N} Y_i^k\big)\ \ge \ \gamma N\big)
\Big).
\end{split}
\ee

\paragraph{Step 3: Dealing with one site.}
We show that for a function $\epsilon(\gamma)$ going
to 0 with $\gamma$, for all $ i\in\Lambda_N$
\be{devil-5}
P\Big(\sum_{k=1}^{t_N} Y_i^k\ge\gamma\log(N)-H\Big)
\ge \exp\Big(-\epsilon(\gamma) \log N\Big).
\ee
We define the successive wave numbers at which the column at 1 grows.
Let $\tau$ be the {\it number of waves}
needed so as to increase by at least one the height of site 1. Then, 
let $\tau_1=\tau$ and 
$\tau_n=\tau\circ\theta(\tau_1+\dots+\tau_{n-1})$. 
Note that for any integer $n$
\be{def-geo1}
P(\tau_1>n) = 
P\big(\tilde Y_i^{(1)}=0,\dots,\tilde Y_i^{(n)}=0\big)
= \exp\big(-n\frac{\KVE H^2x_N}{N}\big),
\ee
where we used \eqref{pikappa}.
Note that at the number of waves $t=\tau_1+\dots+\tau_{k-1}$, 
the configuration
$\tilde\sigma_i^{(t)}$ is larger or equal than $H+k-1$. We have,
using Lemma~\ref{lem-aux}
\be{def-geok}
P\big(\tau_k>n\big| \kG_{k-1}\big)
\le \exp\big(-n\frac{\KVE(H+k-1)^2x_N}{N}\big).
\ee
Then, we are in the setting of Lemma~\ref{lem-comp-geo}, and
have a comparison with independent geometric random variables
$\{\tilde \tau_k,\ k\ge 1\}$ with 
\[
E[\tilde \tau_k]=\frac{N}{\KVE(H+k-1)^2x_N}.
\]
Then, (with the abuse of taking $\gamma\log(N)$ to be integer)
\be{def-geo2}
\begin{split}
P\Big(\sum_{k=1}^{t_N} Y_i^k\ge& \gamma \log(N)-H\Big)\ge 
P\big(\tau_1+\dots+\tau_{\gamma \log N}\le t_N\big)
\ge  P\big(\tilde\tau_1+\dots+\tilde\tau_{\gamma \log N}\le t_N\big)\\
\ge & \prod_{k=1}^{\gamma\log N} 
P\big(\tilde\tau_k\le \frac{t_N}{\gamma \log N} \big)
\ge  \prod_{k=1}^{\gamma\log N} \Big(1-
\exp\big(-t_N\frac{\KVE(H+k)^2x_N}{\gamma N\log(N)}\big)\Big)\\
\ge & \prod_{k=1}^{\gamma\log N}
\Big(1-\exp\big(-\frac{\gamma b\KVE}{2e}
\frac{k^2}{\gamma^2 \log(N)}\big)\Big).
\end{split}
\ee
Now use the estimate $1-e^{-x}\ge e^{-a}x$ for 
$0\le x\le a$ with $a=\gamma b\KVE/2e$. Now, 
\be{geo-3}
\begin{split}
 \prod_{k=1}^{\gamma\log N}
\Big(1-\exp\big(-\frac{\gamma b\KVE}{2e}
\frac{k^2}{\gamma^2 \log(N)}\big)\Big)
\ge & (ae^{-a})^{\gamma log(N)}\frac{(\gamma \log(N))!^2}
{(\gamma^2\log^2(N))^{\gamma \log(N)}}\\
\ge & \exp(-\epsilon(\gamma) \log(N)),
\end{split}
\ee
where we took $\epsilon(\gamma)=\gamma\log(ae^{-a})-2\gamma$.

The proof of the Proposition~\ref{prop-init} is completed if 
$\a>\epsilon(\gamma)$.
\qed
\section{Growing Columns}
\label{sec-column}
In this section we present a simple way to bound the height of
the maximal pile, based on stochastic domination, 
for both ballistic and diffusive deposition.
Indeed, for both models we construct a
sequence of \emph{inter-arrival} times of explorers on a given
column, say column  number $1$, stochastically dominated by
independent geometric variables.

We now state three propositions that are crucial 
in the proofs of Theorems~\ref{main-theo-bal} and \ref{main-theo}.  
The propositions bound the probabilities of
building a high pile, and are proven at the end of this section.

\begin{proposition}
\label{t:lem-ball} 
Consider ballistic deposition. There is a constant $\kappa>0$ such that
for all $H>2$, and any site $i\in G_N$,
\be{conj-1}
P\big(\sigma_i(N)> H\big)\le \exp(-\kappa H).
\ee
\end{proposition}
Proposition~\ref{t:lem-ball} implies that, for a given column,  
$N$ particles are not enough to reach a maximal 
height of order $\log(N)/\kappa$.

\begin{proposition}
\label{t:lem-diff} 
Consider diffusive deposition. For $H<N^{\frac{1}{2}-\epsilon}$,
and any site $i\in G_N$, and $X$ a positive  integer
\be{conj-1-diff}
P\big(\sigma_i(X)> H\big) \le
\exp\Big(\kappa\frac{X}{N}H^2-\frac{\pi}{4}H\Big),
\ee
with $\kappa=1+O(N^{-1/2})$.
\end{proposition}
Proposition~\ref{t:lem-diff} implies that for $\alpha$ small,
$\alpha N/\log(N)$ particles are not enough to reach a maximal 
height of order $\log(N)$.

Finally, we consider diffusive deposition,
with a configuration in $\ER(A)$, and with one distinguished
site, say $i$, above height $\gamma \log(N)$.
We show that as long as we do not leave the {\it early
regime},  see equation (\ref{ear-reggg}) we have a {\it fast growth}. Let $\tau_{\ER^c}$ the time at
which you exit the early regime.

\bp{prop-above}
Let $\chi,\gamma,C$ be any positive constants with $\chi<1/2$.
Assume that there is a distinguished
site, say $i^*$, with $\sigma_{i^*}\ge \gamma \log(N)$.
Then, we have
\be{ineq-above}
P_\sigma\Big(\sigma_{i^*}\big(C \frac{N}{\log(N)}\big)
< N^\chi\big| \tau_{\ER^c}>C \frac{N}{\log(N)}\Big) \le
\exp\Big( -\gamma \big(\kappa(A)\frac{\gamma C}{2}-
\frac{\pi}{2}\big)\log(N)\Big).
\ee
\ep
\subsection{Growing a Column in ballistic deposition}
\label{ch_BD}

\smallskip
\par\noindent
\textit{Proof of Proposition~\ref{t:lem-ball}.\/}
By lemma \ref{lem-DB-upp} 
\be{main-2}
P\big(\tau_1>k\big)\ge \big(1-\frac{1}{N}\big)^k
\qquad\text{so that}\qquad E[\tau_1]> N.
\ee
and in general
\be{main-3}
P\big(\tau_i>k\big)\ge \big(1-\frac{i}{N}\big)^k.
\qquad\text{so that}\qquad E[\tau_i]> \frac{N}{i}.
\ee
This implies that \reff{conj-1} is a large deviation event since
\[
\sum_{i=1}^H E[\tau_i]\ge N\log(H).
\]
More precisely by Lemma~\ref{lem-comp-geo} we have
\[
P\big(\tau_1+...+\tau_H<N )\le
P\big(\tilde\tau_1+...+\tilde\tau_H<N  )\]
with $\{\tilde\tau_i,\ i=1,\dots,H\}$ 
independent geometric variables of mean $\frac{N}{i}$.
By the exponential Chebyshev's inequality we get, for every  $\l>0$ 
\[
P\big(\tilde\tau_1+...+\tilde\tau_H<N  \big)\le
e^{\l N}\prod_{i=1}^H E[e^{-\l\tilde\tau_i}]
\]

Note that for a geometric variable $X$ of mean $1/p$
\[
E[\exp(-\lambda X)]=\frac{p}{e^\lambda-(1-p)}.
\]
When $\lambda$ is positive, $\exp(\lambda)-1\ge \lambda$, and we have
\[
E[\exp(-\lambda X)]\le 1-\frac{\lambda}{\lambda+p}.
\]
Now,
\[
\prod_{i=1}^H E[\exp(-\lambda \tilde\tau_i)]\le 
\prod_{i=1}^H \big(1-\frac{\lambda}{\lambda+i/N}\big)\le
\exp(-\lambda \sum_{i=1}^H \frac{1}{\lambda +i/N}).
\]
We choose $\lambda=i_0/N$ so that by the asymptotic of the harmonic series
\[
\lambda \sum_{i=1}^H \frac{1}{\lambda +i/N}=i_0\sum_{i=i_0+1}^{i_0+H}
\frac{1}{i}\ge i_0\log\big(\frac{i_0+H+1}{i_0+1}\big).
\]
So, we obtain
\[
P\big(\tilde\tau_1+\dots+\tilde\tau_H\le N\big)\le 
\exp\Big(-i_0\big[\log\big(\frac{i_0+1+H}{i_0+1}\big)-1\big]\Big).
\]
Now, if we choose $i_0$ to be the integer part of $\alpha H$, for
some constant $\alpha$, then 
\[
P\big(\tilde\tau_1+\dots+\tilde\tau_H\le N\big)\le
 \exp\Big(-\alpha H(1-\frac{1}{\a H})
\big(\log(\frac{1+ \alpha}{\alpha})-1\big)\Big).
\]
If $\a$ is sufficiently small, say $\a=\frac{1}{2}$, then
$\log((1+ \alpha)/\alpha)-1>0$, and
Lemma~\ref{t:lem-ball} is established.
\qed

\subsection{Growing a Column in diffusive deposition}
\label{ch_DD}

\smallskip
\par\noindent
\textit{Proof of Proposition~\ref{t:lem-diff}.\/}
We follow the arguments of the previous proof.
By using Lemma~\ref{lem-DD-upp} for any $i\ge 1$, and integer $n$
\[
P\big(\tau_i>n\big)\ge \big(1-\kappa\frac{i^2}{N}\big)^n,
\qquad\text{with}\quad \kappa=1+O(N^{-1/2}).
\]
By Lemma~\ref{lem-comp-geo}, we have for any $H,X$
\[
P\big(\tau_1+...+\tau_H<X )\le
P\big(\tilde\tau_1+...+\tilde\tau_H<X ),
\]
with $\{\tilde\tau_i,\ i\ge 1\}$ independent geometric variables with
$E[\tilde \tau_i]=N/(\kappa i^2)$.

Then, for every  $\l>0$, by Chebyshev's inequality
\be{diff-2}
P\big(\tau_1+\dots+\tau_H<X\big)\le
\exp\big(\lambda X-\lambda \sum_{i=1}^H \frac{1}
{\lambda+\kappa i^2/N}\big).
\ee
We set $\lambda=\kappa H^2/N$ and we note that 
\be{diff-3}
H^2\sum_{i=1}^{H} \frac{1}{H^2+i^2}
= \sum_{i=1}^{H} \frac{1}{1+(\frac{i}{H})^2}
\ge H\int_0^1 \frac{dx}{1+x^2}=H\frac{\pi}{4}.
\ee
We conclude obtaining \reff{conj-1-diff}.
\qed

\paragraph{Proof of Proposition~\ref{prop-above}.}

Again, let $\{\tau_1,\tau_2,\dots,\}$ be the random number of
explorers linked with growing a column at $i^*$ from $\sigma$
with an initial state with $\sigma_{i^*}=\gamma\log(N)$.
By Corollary \ref{cor-ear} we have for any integer $m<X$
$$
P_\sigma\big(\tau_k>m\big|\tau_1,\dots,\tau_{k-1},\tau_{\ER^c}>X \big)
\le \Big(1-\kappa(A)\frac{(\gamma\log N+k-1)^2}{N}\Big)^m
=P(\tilde\tau_k>m)
$$
Therefore, by Lemma~\ref{lem-comp-geo}
$$
P_\sigma\big(\sum_{k=1}^{N^\chi}\tau_k>X\big| \tau_{\ER^c}>X \big)\le
P\big(\sum_{k=1}^{N^\chi}\tilde \tau_k>X\big).
$$
By Chebyshev inequality,
$$
\E e^{\lambda\tilde\tau_k}=\frac{p_k}{e^{-\lambda}-(1-p_k)}
=1+\frac{a}{p_k-a},\quad\text{with}\quad
p_k:=\kappa(A)\frac{(\gamma\log N+k-1)^2}{N},
$$
assuming $a:=1-e^{-\lambda}<p_k$. Thus,
$$
\E e^{\lambda\tilde\tau_k}<\exp\{a/(p_k-a)\}
\;\;.
$$
Hence, 
$$
P\Big(\sum_{k=1}^{N^\chi}\tilde \tau_k>X\Big)
\le 
\exp\Big(-\lambda X+\sum_{k=1}^{N^\chi}\frac{a}{p_k-a}\Big)
\;\;.
$$
We choose $a=\kappa(A)(\gamma^2\log^2N)/N$ (and $a\ge \lambda/2$)
and $X=CN/\log N$, and we have
$$
\sum_{k=1}^{N^\chi}\frac{a}{p_k-a}
\le 
\gamma\log N\int_0^\infty\frac{dx}{1+x^2}=\frac{\pi}{2}\gamma\log N,
$$
so that
$$
\begin{array}{rcl}
{\displaystyle
 P\Big(\sum_{k=1}^{N^\chi}\tilde \tau_k>X\Big)
}
&\!\!\le&\!\! 
{\displaystyle
\exp\Big(-\kappa(A)\frac{\gamma^2\log^2N}{2N} 
C\frac{N}{\log N}+\frac{\pi}{2}\gamma\log N\Big)
}\\
&\!\!=&\!\!
{\displaystyle
\exp\Big(-\gamma\log N[\kappa(A)\frac{\gamma C}{2}-\frac{\pi}{2}]\Big)
}
\end{array}
$$
\qed
\subsection{Growing a Tower in diffusive deposition}
\label{s:zoccolo}
In this section, we bound the probability of 
forming a {\it high tower} of explorers in 
diffusive deposition.

Fix a region $\Lambda\subset\{1,\dots,N\}$
of size $L$. Let $H$ be a fixed positive integer,
and $\C=\Lambda\times\{0,\dots,H\}$.
Given $\sigma\in\N^{G_N}$ we define
\begin{equation}
\label{csigma}
\sigma\wedge \C:=\{i\in \Lambda:\sigma_i<H\}\subset G_N,
\quad\text{and note that}\quad
|\sigma\wedge \C|= \sum_{i\in \Lambda} \ind_{\{\sigma_i<H\}} .
\end{equation}
\begin{proposition}
\label{t:zoccolo}
Consider diffusive deposition. 
For any positive $H, X$ and $\xi$  with $X<N/2$, we have that 
\begin{equation}
\label{zoccolo}
P(\forall x\in \Lambda,\ \sigma_x(X)\ge H)
\le \exp\Big(LH\big( \xi \frac{HX}{N}-\xi
\log(\frac{1+\xi+1/H^2}{\xi+1/H^2})\big)\Big)
\end{equation}
Moreover, for any positive $a,\gamma$, and
a positive real $\chi$ satisfying $4\chi<a\gamma^2\exp(-2a\gamma)$,
and the choice $H=\gamma \log(N)$, $X=aT_N$ and $L=N^{1-2\chi}$, we 
have
\be{ineq-tower}
P\big(|\{i:\, \sigma_i(aT_N) >\gamma \log N\}|>N^{1-2\chi}
\big)\le \exp\big(-\chi N^{1-2\chi}\big).
\ee
\end{proposition}

\bpr
Lemma \ref{lem-DD-upp} immediately yields
\begin{displaymath}
P\big(\text{ Explorer attaches to } \sigma\wedge \C\ \big|\ \sigma\big)
\le
\kappa_{A}\sum_{x\in  \sigma\wedge \C}
\frac{\sigma_x^2+1}{N}
\end{displaymath}
Note that
\begin{displaymath}
\sum_{x\in  \sigma\wedge \C}
\sigma_x^2+1 \le L+H\sum_{x\in \sigma\wedge \C}\sigma_x.
\end{displaymath}
We define $n_\C(\sigma)=\sum_{x\in \sigma\wedge \C}\sigma_x$,
so that, for $N$ large enough we have that 
\be{grow-UB}
P\big(\text{ Explorer attaches to } \sigma\wedge \C\ \big|\ \sigma\big)
\le \frac{\kappa_{A}}{N}[L+Hn_\C(\sigma)].
\ee
This allows us  to define, as before, 
geometric random variables stochastically
smaller than the number of explorers needed to 
settle one of them in $\C$.
Let $\tau_1$ be the number of explorers needed in order
that one settles in $\C$, 
when we start with the empty configuration.
By induction, when $k-1$ explorers are settled in $\C$, define
$\tau_k$ to be the number of explorers needed to settle the $k$--th
explorer in $\C$, and we do this up to  time $LH$. Then
for any configuration $\sigma$ with $n_\C(\sigma)=k-1$, for any
positive integer $m$
\be{def-time}
P(\tau_k>m\ |\ \sigma)\ge \big(1-\frac{kH+L}{N}\big)^m=
P(\tilde \tau_k>m).
\ee
We invoke again Lemma~\ref{lem-comp-geo} to obtain
\be{stoch-ineq}
\begin{split}
P\big(\tau_1+\dots+\tau_{HL}\le X\big)\le&
P\big(\tilde\tau_1+\dots+\tilde\tau_{HL}\le X\big)
\le e^{\lambda X} \prod_{k=1}^{HL} E[\exp(-\lambda \tilde \tau_k)]\\
\le & \exp\Big(\lambda X-\lambda \sum_{k=1}^{HL} \frac{1}{\lambda+(kH+L)/N}
\Big)\qquad\text{choose}\quad \lambda=\frac{H^2 L}{N} \xi\\
\le & \exp\Big(\xi\frac{H^2L X}{N}-\xi HL\sum_{k=1}^{HL} \frac{1}
{HL \xi +k+L\H}\Big)\\
\le & \exp\Big(LH\big( \xi \frac{HX}{N}-\xi
\log(\frac{1+\xi+1/H^2}{\xi+1/H^2})\big)\Big).
\end{split}
\ee
With the choice 
$H=\gamma\log(N)$, $L=N^{1-2\chi}$, and $X=aT_N$ we get
\be{main-zoccolo}
\begin{split}
P\big(|\{i:\, \sigma_i(aT_N) >&\gamma \log N\}|>N^{1-2\chi}
\big)\le {N\choose  L} \exp\Big(LH\big( \xi \frac{HX}{N}-\xi
\log(\frac{1+\xi+1/H^2}{\xi+1/H^2})\big)\Big)\\
\le & \exp\Big(L\log(N/L)+L H \frac{HX}{N}\xi-LH\xi
\log(\frac{1+\xi+1/H^2}{\xi+1/H^2})\big)\Big)\\
\le & \exp\Big(-L\log(N)
\big(\gamma \xi\log(\frac{1+\xi}{\xi})-2\chi-a\gamma^2\xi\big)\Big).
\end{split}
\ee
First choose $\xi=\exp(-2a\gamma)$, to get
$$
\big(h\xi\log(\frac{1+\xi}{\xi})-2\chi- xh^2\xi\big)
\ge \big(a\gamma^2\xi -2\chi)\big).
$$
Now choose $4\chi< (a\gamma^2)\exp(-2a\gamma)$, and the inequality
\reff{ineq-tower} is obtained.
\epr

\section{Proof of Theorems \ref{main-theo} and \ref{main-theo-bal}}
\label{sec-proof}

\textit{Proof of Theorem~\ref{main-theo}.\/}
Proof of \eqref{result-1}:
the statement follows immediately by Proposition~\ref{t:lem-diff} 
with $H=3\log(N)$ and $X\le c N/\log(N)$, for $c$ small enough.
Indeed by Proposition~\ref{t:lem-diff}
we have for the complementary event
\be{bbe-1}
\begin{split}
P(\exists i:\;  \sigma_i(X)>3\log N)\le &
\sum_{i=1}^N P( \sigma_i(X)>3\log N)\\
\le &
N \exp\Big(3\log N\big(3\kappa c-\frac{\pi}{4}\big) \Big).
\end{split}
\ee
Hence
\[
P(\exists i:\;  \sigma_i(X)>3\log N) \le N^{1-3(\pi/4-3\kappa c)}.
\]
This concludes the proof
since the exponent of $N$ is less than $-1$ when $9c\kappa<(3\pi/4-2)$.
(recall that $\kappa=1+O(N^{-1/2})$).

\paragraph{Proof of \eqref{result-2}.} 
Recall that from Proposition~\ref{prop-init} there is $b>0$ (and 
\reff{goal-sec5} for the quantitative estimate), so that very likely
$\tau_{\VE^c}<b T_N$, where $T_N=N/\log(N)$. We therefore
condition on the evolution up to $\tau_{\VE^c}$.
\be{glob-1}
\begin{split}
P\big(\max\sigma_x(aT_N+bT_N)<N^\chi\big)\le &
P\big( \tau_{\VE^c}<bT_N,\max\sigma_x(aT_N+bT_N)<N^\chi\big)\\
&\quad+ P\big( \tau_{\VE^c}\ge bT_N\big).
\end{split}
\ee
Now, in the first term use  Markov's property
at time $\tau_{\VE^c}$, calling
for simplicity $\sigma(\ER):=\sigma(\tau_{\VE^c})$
\be{glob-2}
P\big( \tau_{\VE^c}<bT_N,\max\sigma_x(aT_N+bT_N)<N^\chi\big)\\
\le E\Big[\ind_{\{\tau_{\VE^c}<bT_N\}}
P_{\sigma(\ER)}\Big(\max\sigma_x(aT_N)<N^\chi\Big)\Big]
\ee
Now, $\sigma(\ER)\in \ER(A)$, for $A>2\gamma\beta$, and
there is $i^*\in G_N$ such that $\sigma_{i^*}(\ER)\ge \gamma \log(N)$.
Thus,
\be{glob-3}
\begin{split}
P_{\sigma(\ER)}\big(\max_x \sigma_x(aT_N)<N^\chi\big)\le &
P_{\sigma(\ER)}\big(\max_x \sigma_x(aT_N)<N^\chi\big|
\ \tau_{\ER^c}>aT_N\big) P_{\sigma(\ER)}\big( \tau_{\ER^c}>aT_N\big)\\
&\qquad
+P_{\sigma(\ER)}\big(\tau_{\ER^c}<aT_N,\ \max_x\sigma_x(aT_N)<N^\chi\big).
\end{split}
\ee
The first term on the right hand side of \reff{glob-3}
is dealt with by Proposition~\ref{prop-above}.

The next lemma deals with the second term on the right hand side of 
\reff{glob-3}. 
\begin{lemma}
\label{t:esco}
Let $\sigma\in \ER(A)$ be a configuration such that $\max\sigma_x=\gamma\log(N)$,
for some positive $\gamma$. For $a>0$ such that $a\gamma\le A-1$, and 
$4\chi<a\gamma^2\exp(-2a\gamma)$, we have
\begin{equation}
\label{esco}
P_\sigma\big(\tau_{\ER^c}<aT_N,\ \max_x\sigma_x(aT_N)<N^\chi\big)
\le \exp(-\chi N^{1-2\chi}).
\end{equation}
\end{lemma}
\bigskip
\par\noindent
\textit{Proof of Lemma~\ref{t:esco}.\/}
Note that when $a\gamma\le A-1$,
\be{tower-event}
\{\tau_{\ER^c}< aT_N,\ \sigma_x(aT_N)<N^\chi\}\subset
\{\big|\{x:\ \sigma_x(aT_N)>\gamma \log(N)\}\big|>N^{1-2\chi}\}.
\ee
Indeed, on the event $\{\tau_{\ER^c}< aT_N,\ \sigma_x(aT_N)<N^\chi\}$,
\be{esco2}
\begin{split}
A N\le \sum_i\sigma^2_i(\tau_{\ER^c})=&
\sum_{i:\, \sigma_i(\tau_{\ER^c})\le \gamma \log N}
\sigma^2_i(\tau_{\ER^c})+
\sum_{i:\, \sigma_i(\tau_{\ER^c})>\gamma \log N}
\sigma^2_i(\tau_{\ER^c})\\
\le & \gamma (\log N)aT_N+ N^{2\chi} 
|\{i:\, \sigma_i(\tau_{\ER^c}) >\gamma \log N\}|\\
\le & a\gamma N+ N^{2\chi} 
|\{i:\, \sigma_i(\tau_{\ER^c}) >\gamma \log N\}|.
\end{split}
\ee
Proposition~\ref{t:zoccolo} deals with growing large towers,
and by using inequality \reff{ineq-tower} the proof is complete.
\qed

\bigskip
\par\noindent
\textit{Proof of Theorem~\ref{main-theo-bal}.\/}
The statement follows immediately by Proposition~\ref{t:lem-ball} 
with $H=\gamma\log(N)$ with $\gamma>2/\kappa$.
Indeed, by Proposition~\ref{t:lem-ball}, there exists a constant $\kappa$ 
such that
\[
P(\exists i:\;  \sigma_i(N)>\gamma\log N)\le 
N e^{-\kappa \gamma \log N} =N^{1-\kappa \gamma},
\]
and the proof concludes.
\qed

\noindent{\bf Acknowledgements.}
A.A. Thanks Robin Pemantle for discussions on urns.
A.A. and E.S. thank the CIRM for a friendly atmosphere
during their stay as part of a {\it research in pairs} program.
This work has been carried out thanks to the support of
A$^*$MIDEX grant
(ANR-11-IDEX-0001-02) funded by the French Government
"Investissements d'Avenir" program.
E.S. thanks Universit\'e Paris--Est, Cr\'eteil.
Finally, we thank two anonymous referees 
for their careful reviewing.

\newpage

A. Asselah\\
Aix-Marseille Universit\'e, Marseille, \\
and LAMA, Universit\'e Paris-Est Cr\'eteil, France.\\

E.N.M. Cirillo, \\
Dipartimento di Scienze di Base e Applicate per l'Ingegneria, \\
Sapienza Universit\`a di Roma, Roma, Italy. \\

E. Scoppola, \\
Dipartimento di Matematica e Fisica, \\
Universit\`a di Roma 3, Roma, Italy.\\

B. Scoppola, \\
Dipartimento di Matematica, \\
Universit\`a di Tor Vergata, Roma, Italy.\\

\end{document}